\documentclass[12pt]{amsart}

\usepackage{latexsym}
\usepackage[all]{xy}
\usepackage{amsmath, amsthm }
\usepackage{amssymb}
\usepackage{amsfonts}
\usepackage{color}
\usepackage{mathtools}

\usepackage{hyperref}
\usepackage{xcolor}
\hypersetup{colorlinks=true,urlcolor=blue,linkcolor=blue,citecolor=blue}

\usepackage{cite}

\newcommand {\Perf} {\mathsf{Perf}}
\newcommand {\QCoh} {\mathsf{QCoh}}
\newcommand {\Coh} {\mathsf{Coh}}
\newcommand {\Vect} {\mathsf{Vect}}

\newcommand {\Map} {\mathbf{Map}}
\newcommand {\Parf} {\mathsf{Perf}}
\newcommand {\rh} {\mathbb{R}\underline{Hom}}
\newcommand {\rch} {\mathbb{R}\underline{\mathcal{H}om}}
\newcommand {\OO} {\mathcal{O}}
\newcommand {\DR}{\mathbf{DR}}
\newcommand {\Fol}{\mathcal{F}ol}

\newcommand {\F} {\mathcal{F}}
\newcommand {\G} {\mathcal{G}}
\newcommand {\A} {\mathcal{A}}

\newcommand {\T} {\mathbb{T}}
\newcommand{\C}{\mathbb{C}}
\newcommand{\Z}{\mathbb{Z}}
\newcommand{\LL}{\mathbb{L}}

\newcommand {\Spf} {\mathsf{Spf}}
\newcommand  {\dgart}     {\mathbf{dgart}}
\newcommand  {\dg}     {\mathbf{dg}}

\newcommand  {\medg}     {\epsilon-\mathbf{dg}^{gr}}
\newcommand {\scat} {\mathbf{Cat}_{\infty}}
\newcommand  {\mecdga}     {\epsilon-\mathbf{cdga}^{gr}}

\newcommand  {\cdga}     {\mathbf{cdga}}

\newcommand  {\dAff}     {\mathbf{dAff}}
\newcommand  {\dAfff}     {\mathbf{dAff}^\mathcal{F}}

\newcommand  {\dSt}   {\mathbf{dSt}}

\newcommand{\s}{\infty}
\newcommand{\HH}{\mathbb{H}}

\newcommand{\D}{\mathcal{D}}

\newcommand{\fcirc}{\widehat{S}^1}

\theoremstyle{plain}
\newtheorem{thm}{Theorem}[subsection]
\newtheorem{df}[thm]{Definition}
\newtheorem{prop}[thm]{Propositon}
\newtheorem{rmk}[thm]{Remark}
\newtheorem{cor}[thm]{Corollary}
\newtheorem{ex}[thm]{Example}
\newtheorem{lem}[thm]{Lemma}

\title[\resizebox{4.5in}{!}{Algebraic foliations and derived geometry I: 
The Riemann-Hilbert correspondence}]{Algebraic foliations and derived geometry I: 
The Riemann-Hilbert correspondence}

\author{Bertrand To\"en and Gabriele Vezzosi}
\date{May 2020}

\begin{document}

\maketitle

\begin{abstract} 
This is the first in a series of papers about foliations in derived geometry. After introducing derived foliations on arbitrary derived stacks, we concentrate on quasi-smooth and rigid derived foliations on smooth complex algebraic varieties and on their associated formal and analytic versions. Their truncations are classical singular foliations defined in terms of differential ideals in the algebra of forms. We prove that a quasi-smooth rigid derived foliation on a smooth complex variety $X$ is formally integrable at any point, and, if we suppose that its singular locus has codimension $\geq 2$, its analytification is a locally integrable singular foliation on the associated complex manifold $X^h$. We then introduce the derived category of perfect crystals on a quasi-smooth rigid derived foliation on $X$, and prove a Riemann-Hilbert correspondence for them when $X$ is proper. We discuss several examples and applications.
\end{abstract}

\tableofcontents

\section*{Introduction}

This is the first of a series of works on foliations (mainly algebraic and holomorphic)
and derived geometry. In this paper we present a notion of a \emph{derived foliation}
on algebraic or holomorphic varieties, that we think is interesting for the study of 
foliations with singularities. The point of view adopted here is not completely new 
and goes back to previous 
works by Tony Pantev and the authors on existence of potentials for
shifted symplectic structures (see e.g. \cite{tony}). In a nutshell,  
a \emph{derived foliation} $\F$ on a scheme $X$ consists of a perfect complex $\LL_{\F}$ on $X$ together
with a map $a : \OO_X \longrightarrow \LL_{\F}$ that satisfies formal properties of being
a de Rham differential (i.e. is a derivation squaring to zero). One major
difficulty is to define the precise higher coherences for such a structure, encoding the
fact that $a^2$ does not really identically vanish but is only homotopic to zero in a homotopy 
coherent way. This is achieved by defining derived foliation as \emph{graded mixed 
commutative differential graded algebras} (graded mixed cdga's, for short) satisfying some extra properties (see Definition 
\ref{d1}).

In this work we quickly restrict to the case of \emph{quasi-smooth derived foliations} $\F$,
which consists of restricting $\LL_{\F}$ to be just a two terms complex
of vector bundles. Among derived foliations, these quasi-smooth derived foliations  are the closest to classical foliations in the usual sense, 
and we think they form the most important class of derived foliations. 
A quasi-smooth derived foliation $\F$ on a smooth variety $X$ can be \emph{truncated} into a usual algebraic 
\emph{singular foliation} $\tau_0(\F)$ 
on $X$ (e.g. in the sense of \cite{ba, ay}). More precisely, the kernel of the morphism 
$\Omega_X^1 \longrightarrow H^0(\LL_\F)$, induced by $a$, defines a differential 
ideal  inside differential forms and thus a singular foliation $\tau_0(\F)$ on $X$. We remark however that arbitrary 
singular foliations are not derived foliations: they can be represented by graded mixed
algebras, but these do not satisfy our conditions (except if the foliation has no singularities).
Notice also that being the truncation 
of a derived foliation is a non-trivial condition, even locally in the analytic topology. 
Therefore, derived foliations are not really generalizations of singular foliations, and these two class
of objects do not live in the same categories. Rather, it is more useful to keep in mind the intuition that 
derived foliations are additional structures on their truncated singular foliations
making them better behaved objects.\\

The first two main results of this work are the following integrability theorems. Note that 
for a derived quasi-smooth foliation being integrable, i.e. being 
induced by a morphism between smooth varieties, implies that its truncated singular foliation is
also integrable (by the same morphism). However, the converse is in general wrong.

\begin{thm}[Prop. \ref{p6} and Cor. \ref{cp6'}]\label{ti1}
Let $X$ be a smooth variety and $\F$ be a quasi-smooth derived foliation on $X$.
Assume that $\F$ is rigid (i.e. the induced map $H^0(a): \Omega_X^1 \rightarrow H^0(\LL_{\mathcal{F}})$ is surjective). Then
\begin{enumerate}
\item The derived foliation $\F$ is formally integrable around each point $x \in X$.

\item If we further assume that $\F$ has no singularities outside a closed subset
of codimension at least $2$, then $\F$ is analytically integrable, locally in the analytic topology on $X$.
As a consequence
the truncated singular foliation $\tau_0(\F)$
is analytically integrable, locally in the analytic topology on $X$.
\end{enumerate}
\end{thm}

Part $(1)$ of the above theorem is a consequence (Cor. \ref{cp4}) of a more general result concerning the \emph{local 
structure} of quasi-smooth derived foliations (see Proposition \ref{p4}), while part $(2)$ is a consequence
of $(1)$ and of a theorem of Malgrange (\cite{mal2}). We remark here that a consequence of the
above result is that a singular foliation
which is not formally integrable locally at all points \emph{can not} be the truncation of a derived quasi-smooth foliation. \\

The second main result of this work is a \emph{Riemann-Hilbert correspondence for derived quasi-smooth
foliations}. We first introduce the notion of a \emph{crystal} along a derived foliation $\F$, which 
morally consists of a vector bundle together with a partial connection along the leaves
of $\F$. Once again,
there are homotopical coherences to be taken into account, and crystals are rather defined
as certain graded mixed dg-modules over the graded mixed dg-algebra defining the derived foliation.
On the other hand, a derived foliation $\F$ defines a sheaf $\OO_{\F^h}$, in the analytic topology,
of locally constant functions along $\F$. This is a sheaf of commutative dg-algebras, which is in general
not concentrated in degree zero, and whose higher cohomology sheaves reflect the singularities
of $\F$. The Riemann-Hilbert correspondence can then be stated as follows (see Cor. \ref{ct1}):

\begin{thm}\label{ti2}
Let $\F$ be a quasi-smooth and rigid foliation on a smooth and proper algebraic variety $X$.
Assume that $\F$ is non-singular outside of a closed
subset of codimension at least $2$. There is an equivalence of categories
$$\Vect(\F) \simeq \Vect(\OO_{\F^h})$$
between on the l.h.s. the category of crystals along $\F$, and on the r.h.s.
the category of sheaves of $\OO_{\F^h}$-dg-modules which are locally free of finite rank.
\end{thm}

The above theorem is a consequence of two results proved in the text: a more general statement 
(valid without the rigidity or
codimension assumptions) which relates perfect complexes of crystals with a nilpotent condition 
and perfect complexes of $\OO_{\F^h}$-dg-modules (see Theorem \ref{t1}), and the fact that a vector 
bundle crystal (i.e. an object in $\Vect(\F)$) is nilpotent once $\F$ satisfies the hypotheses of 
Theorem \ref{ti2} (see Theorem \ref{t0}). We also prove that
the above theorem is compatible with cohomologies, giving rise to an isomorphism between 
algebraic de Rham cohomology  along the leaves of a crystal and the analytic cohomology of the corresponding
sheaf of $\OO_{\F^h}$-module.
Note that Theorem \ref{ti2} recovers for instance Deligne's relative Riemann-Hilbert correspondence (see \cite{del})
and extends it to the case of a possibly non-smooth morphism (see Section \ref{deligne} for details). 
It is also possible to recover from Theorem \ref{t1} Kato-Nakayama's logarithmic Riemann-Hilbert correspondence (see Section \ref{log}).
\bigskip

\textbf{Plan of the paper.} The present work is organized in four parts. In the first section we present the 
notion of derived foliations. We present several examples and show a formal structure theorem for
quasi-smooth and rigid derived foliations. We also discuss the notion of leaves in the setting. The second section
is devoted to the analytic aspects of derived foliation. We construct the analytification functor and
discuss local integrability, in the analytic topology, of derived foliations. The third section contains the 
definition of crystals along a derived foliations, their analytifications, as well as the notion of nilpotent crystals.  Finally,
the last section contains the statement and proof of the Riemann-Hilbert correspondence. We also have
included some examples and applications. \\

\textbf{Related works.} In \cite{bsy}, the authors borrow their definition of derived foliation from \cite{tony}, and study Lagrangian derived foliations in relation with the problem of realizing the moduli space of sheaves on a Calab-Yau fourfold as the derived critical locus of a (shifted) potential. In \cite{ay}, J. Ayoub have systematically studied \emph{underived} singular foliation on schemes; his theory lives algebraic geometry rather than in derived geometry, and his purposes are somehow different, being related to differential Galois theory.  \\

\textbf{Acknowledgements.} We thank Tony Pantev for useful discussions about derived foliations over the years.
This project has received funding from the European Research Council 
under the European Union's Horizon 2020 research and innovation programm 
(grant agreement No. 2016-ADG-741501).\\

\textbf{Conventions and notations.}
Everything, like vector spaces, commutative dg-algebras (often shortened as cdga's), algebraic varieties etc., is defined over 
the field $\C$ of complex numbers. The $\s$-category of complexes (of $\C$-vector spaces) 
is denoted by $\dg$, and the $\s$-category of topological spaces by $Top$.

By convention $\dAff$ is the $\infty$-category of derived affine schemes locally of finite presentation
over $\C$. Derived Artin stacks are, by definition,  locally of finite presentation. For 
a derived stack $F$ we denote by $QCoh(F)$ its $\s$-category of quasi-coherent complexes. In the same
manner $\Parf(F) \subset QCoh(F)$ denotes the full sub-$\s$-category of 
perfect complexes on $F$.

All the various functors, $Sym$, $\otimes$, $\wedge$, $f_*$, $f^*$, etc. will be suitably 
derived when necessary. We will occasionally need underived functors for which 
we will use specific notations $Sym^u$, $\otimes^u$, $f_*^u$, etc., if necessary.

A vector bundle on $X$ will be a locally free $\mathcal{O}_X$-Module of finite rank.

\section{Derived Algebraic foliations}

In this section, after some reminders on mixed graded structures, we define derived foliations on 
arbitrary derived stacks, give several classes of examples of derived foliations, study derived 
foliations on formal completions, and finally establish the local structure of quasi-smooth rigid 
derived foliations. We also discuss the notion of formal, algebraic and analytic leaves of
derived foliations in general.

\subsection{Reminders on graded mixed objects}\label{reminder}

We remind from \cite{cptvv} (see also the digest  \cite{pv}) the $\s$-category of graded mixed 
complexes (over $\C$).
Its objects are $\Z$-graded objects $E=\oplus_{n\in \mathbb{Z}} E(n)$, inside the category of
cochain complexes  together
with extra differentials $\epsilon_n : E(n) \longrightarrow E(n+1)[-1]$.
These extra differentials  combine into
a morphism of graded complexes $\epsilon : E \longrightarrow E((1))[-1]$  
(where $E((1))$ is the graded complex obtained from $E$ by shifting the weight-grading
by $+1$), satisfying $\epsilon^2=0$. The datum of $\epsilon$ is called
a \emph{graded mixed structure} on the graded complex $E$. The complex $E(n)$ 
is itself called the \emph{weight n 
part} of $E$. 

Morphisms of
graded mixed complexes are defined in an obvious manner, and among them, the
quasi-isomorphisms are the morphisms inducing quasi-isomorphisms on all 
the weight-graded pieces \emph{individually}. 
By inverting quasi-isomorphisms, 
graded mixed complexes constitute an $\s$-category denoted by $\medg$. Alternatively, 
the $\s$-category $\medg$ can be defined as the $\s$-category of quasi-coherent
complexes $\QCoh(B\mathcal{H})$, over the classifying stack 
$B\mathcal{H}$ for the group stack $B\mathbb{G}_a \rtimes \mathbb{G}_m$
(see \cite[Rmk. 1.1.1]{cptvv} and \cite[Prop. 1.1]{pato}).

The $\s$-category $\medg$ comes equipped with a canonical symmetric monoidal 
structure $\otimes$. It is defined on objects by the usual 
tensor product of $\Z$-graded complexes (taken over the base field $\C$), with the mixed structure
defined by the usual formula $\epsilon \otimes 1 + 1\otimes\epsilon$ (see \cite[\S 1.1]{cptvv}). 
When viewed as $\QCoh(B\mathcal{H})$, this is the usual
symmetric monoidal structure on quasi-coherent complexes on stacks. \\

Commutative algebras in $\medg$ form themselves an $\s$-category $\mecdga$, whose
objects are called graded mixed cdga's. Its objects can be described as
$\Z$-graded cdga's $A=\oplus_n A(n)$, endowed with a graded mixed structure $\epsilon$
which is compatible with the multiplication in $A$ (i.e. is a graded biderivation). The 
fundamental
example of such a graded mixed cdga is given by the \emph{de Rham algebra}. For a cdga $A$ 
we can consider its dg-module $\Omega_A^1$ of dg-derivations as well as its
symmetric cdga $Sym^u_A(\Omega^1_A[1])$.
The usual 
de Rham differential induces a graded mixed structure on $Sym^u_A(\Omega_A^1[1])$ 
making it into a graded mixed cdga for which the induced morphism $\epsilon : A \longrightarrow \Omega_A^1$ 
is the usual universal derivation. Applied to a cofibrant model $A'$ of $A$ we get 
a graded mixed cdga $\DR(A):=Sym^u_{A'}(\Omega_{A'}^1[1])$ which is functorial, in the sense of $\s$-categories, in $A$. 
This defines an $\s$-functor
$$\DR : \cdga \longrightarrow \mecdga$$
which can be checked to be the left adjoint to the forgetful $\s$-functor sending 
a graded mixed cdga $A$ to its weight $0$ part $A(0)$. \\

We remind the existence of the realization $\s$-functor
$$|-| : \medg \longrightarrow \dg$$
given by $\rh(\C,-)$, where the derived hom is taken as graded mixed complexes, 
and $\C$ is equiped with the trivial graded mixed complexes purely
concentrated in weight and cohomological degree $0$. The object $\C$ being the
unit of the symmetric monoidal structure on $\medg$, the $\s$-functor $|-|$ 
possesses a natural lax monoidal structure and thus
sends graded mixed cdga's to cdga's. It can be explicitly described as follows. 
For a graded mixed complex $E$ we from the product
$$|E|:=\prod_{p\geq 0}E(p)[-2p]$$
and endow $|E|$ with the total differential $d+\epsilon$, where
$d$ is the cohomological differential of $E$ and $\epsilon$ is the
graded mixed structure. When $E$ is a graded mixed cdga the
formula above for $|E|$ can also be used to describe the multiplicative structure, 
which is induced by the natural maps $E(n)[-2n]\otimes E(m)[-2m] \longrightarrow E(n+m)[-2n-2m]$. \\

\begin{rmk}
The following simple observations will be useful in the rest of the paper.

\begin{itemize}

\item For $A \in \cdga$, the underlying graded cdga of $\DR(A)$, obtained
by forgetting the mixed structure, is naturally equivalent to $Sym_A(\LL_A[1])$, where
$\LL_A$ is the cotangent complex of $A$.

\item As a consequence of the comment above, when $A$ is a smooth algebra, the graded mixed cdga 
$\DR(A)$ 
is canonically equivalent to the usual de Rham algebra $Sym_A(\Omega_A^1[1])$
endowed with its usual de Rham differential as graded mixed structure.

\end{itemize}
\end{rmk}

The notions of graded mixed complexes, graded mixed cdga's and 
de Rham algebras $\DR$ as defined above, all make sense \emph{internally}  
to a (nice enough) base symmetric monoidal $\C$-linear $\s$-category
(see \cite[Section 1.3.2]{cptvv},  as well as \cite[Rmk 1.5 and Section 2.1]{pv}).
These internal notions and constructions can be understood simply as follows.
Graded mixed cgda's and modules make sense over any derived stack $F$, as
quasi-coherent sheaves of $\OO_F$-linear graded mixed cdga's and modules. Equivalently 
the $\s$-category of graded mixed modules over a derived $F$ can be defined as
$\QCoh(F\times B\mathcal{H})$, where, as above,  $\mathcal{H}$ is the group
stack $B\mathbb{G}_a \rtimes \mathbb{G}_m$. Graded mixed cdga's are then 
naturally defined as commutative ring objects inside the symmetric monoidal 
$\s$-category $\QCoh(F\times B\mathcal{H})$. 

Any commutative ring $A$ in $\QCoh(F)$ will be 
called an $\OO_F$-cdga. Any such
$\OO_F$-cdga
possesses an \emph{internal 
de Rham complex}, which is a graded mixed cdga over $F$. We denote this
object by $\DR^{int}(A)$. Moreover, we can apply the direct image functor 
along $F\times B\mathcal{H} \longrightarrow F$ to get a lax
monoidal $\s$-functor
$\QCoh(F\times B\mathcal{H}) \longrightarrow \QCoh(F).$
This lax monoidal $\s$-functor is called the 
\emph{realization} $\s$-functor and is denoted by
$$|-| : \QCoh(F\times B\mathcal{H}) \longrightarrow \QCoh(F).$$
When $A$ is a cdga over $F$, we have a graded mixed cdga $\DR^{int}(A)$
over $F$, and by applying $|-|$ we get a cdga
denoted by $\DR(A):=|\DR^{int}(A)|$, and called
the \emph{de Rham cohomology of $A$ over $F$}. There is also a
relative version, for a morphism $A \longrightarrow B$ of cdga's over $F$, 
which is $|\DR^{int}(B/A)|$, another cdga. The explicit formula
giving  the realization recalled earlier is also valid in this 
internal setting. Indeed, for an object $E \in \QCoh(F\times B\mathcal{H})$, 
its realization $|E|$ is the object in $\QCoh(F)$ given by 
$$|E|=\prod_{p\geq 0}E(p)[-2p]$$
endowed with the total differential, sum of the cohomological differential
and the de Rham differential.

This discussion applies in particular to $F=B\mathcal{H}$ itself.
We have to note here that $\QCoh(B\mathcal{H} \times B\mathcal{H})$
consists of doubly graded mixed complexes, i.e. complexes
endowed with two extra gradings and two associated graded mixed structures
compatible with each others. By our convention the realization
$$|-| : \QCoh(B\mathcal{H} \times B\mathcal{H}) \longrightarrow \QCoh(B\mathcal{H})$$
consists of realizing the first graded mixed structure. For example, 
if one starts with an algebra $A$ in $\QCoh(B\mathcal{H})$ (i.e.  
a graded mixed cdga), then $|\DR^{int}(A)|$ is another
graded mixed cdga. It is obtained by considering $\DR^{int}(A) \in \QCoh(B\mathcal{H} \times B\mathcal{H})$
and realizing it with respect to the \emph{internal mixed structure}, that is the
one induced from the graded mixed structure on $A$ as opposed to the one given by the
de Rham differential. Using the correct convention here is essential 
for the rest of the paper.

If we have a morphism of graded mixed cdga's $A \longrightarrow B$, 
the above construction produces an internal graded mixed 
cdga $\DR^{int}(B/A)$ inside graded mixed complexes. Its realization is
thus a graded mixed cdga $\DR(B/A)$ called the \emph{
internal de Rham cohomology of $B$ relative to $A$}.  \\

With these notations, we have the following lemma recovering a class of graded mixed cdga's $A$ from their $\DR^{int}(A(0)/A)$.  We will use this lemma very often in the rest of the text.

\begin{lem}\label{l0}
Let $A$ be a graded mixed cdga and assume that the canonical morphism 
$$Sym_{A(0)}(A(1)) \longrightarrow A$$
is a quasi-isomorphism of graded cdga's.
Then, the canonical morphism of graded mixed cdga's
$$A \longrightarrow |\DR^{int}(A(0)/A)|$$
is a quasi-isomorphism.
\end{lem}

\noindent\textbf{Proof.} Let $B=A(0)$ and $E=A(1)$. 
The internal cotangent complex of $B$ relative 
to $Sym_B(E)$ is identified with $E[1]$. 
The internal de Rham algebra $\DR(B/A)$ is then 
equivalent to $Sym_B(E[2])$. We are interested in
realizing the internal graded mixed structure coming from
the one of $A$. As $E$ is pure of weight $1$, the induced
graded mixed structure on $E$ is trivial. The same is true
for $Sym^p(E[2])$, and we thus conclude that the internal
graded mixed structure on $Sym_B(E[2])$ is trivial. Since we are realizing internally, 
we have to realize each graded piece individually. But the 
realization of a graded mixed complex $M$ which is pure of weight $p$ 
 is simply $M[-2p]$. 
Therefore, the
realization of the internal de Rham algebra is tautologically given by
$$|\DR^{int}(B/A)|\simeq \oplus_{p\geq 0}|Sym^p(E[2])| = \oplus_{p\geq 0}Sym^p(E).$$

\hfill $\Box$ \\

\subsection{Derived algebraic foliations as graded mixed cdga's}

In this section we present a very general notion of derived
foliations over general derived stacks. Though later in this paper, we will only 
be dealing with derived foliations over smooth varieties, we have decided to  
give a general definition for further record and applications.

\begin{df}\label{d1}
An \emph{affine derived foliation} is a graded mixed cdga $A$ satisfying the following 
extra conditions.
\begin{enumerate}
\item (Connectivity) The underlying cdga $A(0)$ is cohomologically concentrated in non-positive 
degrees and is finitely presented over $\C$.
\item (Perfectness) The $A(0)$-dg-module $A(1)[-1]$ is perfect and connective.
\item (Quasi-freeness) The natural morphism of graded cdga's
$$Sym_{A(0)}(A(1)) \longrightarrow A$$
is a quasi-isomorphism of graded cdga's.
\end{enumerate}
For a derived foliation $A$ as above, the derived affine scheme $X=Spec\, A(0)$ 
is called the \emph{underlying derived scheme of the foliation}, and we will say that the 
foliation is given \emph{over $X$}. The perfect complex
on $X$ determined by $A(1)[-1]$ is called the \emph{cotangent complex
of the foliation}. 

\end{df}

\begin{ex} \emph{Let $X=Spec \, R$ be a smooth affine $\C$-scheme, $TX$ its tangent bundle, and 
$V\subseteq TX$ a sub-bundle whose local sections are closed under the Lie bracket canonically 
defined on local vector fields (i.e. on local sections of $TX/X$). It is well known that if 
$\mathcal{V}^{\vee}$ denotes the $R$-module of local sections of the dual vector bundle 
$V^{\vee}$, then the Lie bracket on local sections of $V$ induces a differential on 
$A:=Sym_R(\mathcal{V}^{\vee}[1])$. This gives $A$ the structure of a derived foliation over $X$.
Therefore, an algebraic foliation in the usual sense can be seen as a derived foliation.}
\end{ex}

More general examples of derived foliations will be given later in this Section.\\

Affine derived foliations form an $\s$-category as follows. Consider the
$\s$-category $(\mecdga)^{op}$, opposite to the $\s$-category of graded mixed cdga's. 
The $\s$-category of affine derived foliations is defined to be the full sub-$\s$-category
of $(\mecdga)^{op}$ consisting of the graded mixed cdga's satisfying the conditions of 
definition \ref{d1}. This $\s$-category will be denoted by $\dAfff$.

We have a canonical $\s$-functor
$$\dAfff \longrightarrow \dAff$$
sending an affine derived foliation $A$ to the derived affine scheme $Spec\, (A(0))$. 

\begin{prop}\label{p1}
The above $\s$-functor
is fibered in the sense of \cite[\S 2.3]{chern}. Moreover, the corresponding $\s$-functor
$$\Fol : \dAff^{op} \longrightarrow \scat$$
is a stack for the \'etale topology. 
\end{prop}

\noindent\textbf{Proof.}
By construction, the $\s$-category $\Fol(Spec\, A)$ is equivalent to 
the opposite $\s$-category of graded mixed cdga's $C$ satisfying the
conditions of Definition \ref{d1} and equipped with a cdga quasi-isomorphism 
$C(0)\simeq A$. The $\s$-category has two distinguished objects, the final 
and initial objects. The final object is $A$ itself, considered
as a graded mixed cdga's purely in weight $0$ with zero graded mixed structure. 
On the other hand, the initial object is $\DR(A)$. 

Let now $f: X=Spec\, A \longrightarrow Y=Spec\, B$
be a morphism of derived affine schemes corresponding to a morphism
of cdga's $B \longrightarrow A$. The pull-back $\s$-functor
$$f^* : \Fol(Y) \longrightarrow \Fol(X)$$
can be understood as follows. Let $\mathcal{F} \in \Fol(Y)$ be 
an object corresponding to a graded mixed cdga $C$ satisfying 
the conditions of Definition \ref{d1} and equipped with a
quasi-isomorphism $C(0) \simeq B$. Associated to $\mathcal{F}$ is a natural diagram of graded
mixed cdgas
$$\xymatrix{
\DR(B) \ar[r] \ar[d] & C \\
\DR(A). & }$$
The pull-back foliation $f^*(\mathcal{F}) \in \Fol(X)$ is then 
given by the graded mixed cdga $C\otimes_{\DR(B)}\DR(A)$. 
This indeed satisfies the conditions of Definition \ref{d1} since it is equivalent,
as a graded cdga, to $Sym_{A}(E)$, where 
$E$ is the following push-out in $B$-dg-modules
$$\xymatrix{B(1) \ar[r] \ar[d] & C(1) \ar[d] \\ A(1) \ar[r] & E}$$ where $A(1)$ is viewed as a 
$B$-dg-module via the map $B \to A$.
This proves the first statement in the proposition, and moreover provides an explicit 
description of pull-back $\s$-functors. This description in turns easily 
implies that the $\s$-functor $\Fol$ is a stack for the \'etale topology, as this
reduces to the fact that quasi-coherent modules is a stack for the \'etale topology.
\hfill $\Box$ \\

The above proposition can be used, by Kan extension along $\dAff^{op} \to \dSt^{op}$, in order to define derived foliations over any base derived stack.

\begin{df}\label{d2}
Let $X \in \dSt$. The $\s$-category 
$$\Fol(X):= \lim_{Spec\, A \to F}\Fol(Spec\,A)$$ is called the \emph{$\s$-category of derived foliations over $X$.}
\end{df}

We note here that when $X$ is a derived DM-stack, then $\Fol(X)$
can be described as the limit $\lim_{U}\Fol(U)$, where $U$ runs
over all derived affine schemes \emph{\'etale over} $X$. By the explicit 
description of pull-backs given in the proof of  Proposition \ref{p1}, 
we see that an object in this limit can be simply represented by a
sheaf of graded mixed cdga's $\A$ over the small \'etale site $X_{\textrm{\'et}}$ of $X$, 
together with an equivalence $\A(0)\simeq \OO_X$, and
satisfying the following two conditions.

\begin{itemize}
\item The sheaf of $\OO_X$-dg-modules $\A(1)[-1]$ is perfect
and connective.

\item The natural morphism of sheaves of graded cdgas 
$$Sym_{\OO_X}(\A(1)) \longrightarrow \A$$
is a quasi-isomorphism.
\end{itemize}

This simple description in terms of sheaves of cdga's is not valid anymore 
for derived foliation over more general derived stacks, such as 
derived Artin stacks for instance. We will quickly restrict ourselves 
to derived foliations over smooth. \\

We introduce the following notations. 

\begin{df}\label{d3}
Let $X$ be a
derived DM-stack and $\F \in \Fol(X)$ be a derived foliation over $X$.

\begin{itemize}
\item The sheaf of graded mixed cdga's $\A$ over $X$ corresponding to 
$\F$ is called the \emph{de Rham algebra along $\F$}. It is denoted by 
$\DR(\F)$.

\item The perfect complex $\DR(\F)(1)[-1]$ over $X$ is called 
the \emph{cotangent complex of $\F$} and is denoted by $\LL_{\F}$.
We thus have a quasi-isomorphism of quasi-coherent sheaves of graded cdga's
over $X$
$$\DR(\F)\simeq Sym_{\OO_X}(\LL_\F[1]).$$
\end{itemize}
\end{df}

Before giving some examples of derived foliations, we fix the following terminology.

\begin{df}\label{d4}
Let $X$ be a
derived DM-stack, $\F \in \Fol(X)$ be a derived foliation over $X$
and $\mathbb{L}_\F \in \QCoh(X)$ is cotangent complex.

\begin{itemize}

\item We say that the foliation $\F$ is \emph{smooth}
if $\LL_\F$ is quasi-isomorphic to a vector bundle on $X$ sitting in degree $0$.

\item We say that the foliation $\F$ is \emph{quasi-smooth}
if $\LL_\F$ is quasi-isomorphic to perfect complex of amplitude contained in 
cohomological degrees $[-1,0]$.

\item We say that the foliation $\F$ is 
\emph{rigid} if the induced morphism of coherent sheaves
$$H^0(\LL_X) \longrightarrow H^0(\LL_\F)$$
is surjective.
\end{itemize}
\end{df}

\begin{rmk} \emph{Definition \ref{d2} above can be extended to more general settings. To start with, we may allow $X$ being any derived Artin stack, and we may furthermore drop 
the connectivity assumption on $\LL_F$ in order to define
\emph{non-connective} derived foliations. These are useful for instance in the setting of shifted symplectic and
Poisson structures (see e.g. \cite{tony,bsy}), but will not be considered in the present work.}
\end{rmk}

To finish this section, we describe a more geometric interpretation of derived foliations in terms
of derived loop spaces and their natural circle action. This point of view will not be used
further in the present paper, but has the advantage of being useful in some contexts, and comes handy. Moreover, using 
the graded circle of \cite{mrt} instead of the formal circle, makes it possible to 
extend the notions of derived foliations outside of the characteristic zero context.

We consider the formal additive group $\widehat{\mathbb{G}}_a$, as well
as its classifying stack $\fcirc:=B\widehat{\mathbb{G}}_a \in \dSt$. The group
$\mathbb{G}_m$ acts on the formal group $\widehat{\mathbb{G}}_a$
and thus on the formal circle $\fcirc$. The stack $\fcirc$ is itself
a group stack and thus acts on itself by translation. These two actions combine
into an action of the group stack $\mathcal{H}:=\fcirc \rtimes \mathbb{G}_m$.
As explained in \cite[Prop. 1.3]{pato} the symmetric monoidal $\s$-category $\QCoh(B\mathcal{H})$
is naturally equivalent to the symmetric monoidal $\s$-category of graded
mixed complexes. 

For a derived DM-stack $X \in \dSt$, its formal derived loop stack 
is defined by
$$\mathcal{L}^fX :=\Map(\fcirc,X).$$
It comes equipped with a canonical action of $\mathcal{H}$. By 
the equivalence recalled above, between $\QCoh(B\mathcal{H})$ and
graded mixed complexes, we see that a derived foliation 
over $X$ is the exact same thing as a 
a derived stack $\F$ over $\mathcal{L}^fX$, together with 
an $\mathcal{H}$-action on $\F$ covering the canonical action 
on $\mathcal{L}^fX$ and such that $\F$ is relatively affine
over $X$ and of the form $Spec_{\OO_X}(\mathbb{L}_{\F}[1])$
(compatibly with the grading where $\mathbb{L}_{\F}$
is of weight one)
for $\mathbb{L}_\F$ a connective perfect complex over $X$. As a result, 
$\Fol(X)$ can be realized as a full sub-$\s$-category 
of $(\dSt/\mathbb{L}^fX)^\mathcal{H}$, of $\mathcal{H}$-equivariant
derived stacks over $\mathbb{L}^fX$. 

The above interpretation of derived foliations makes \emph{pull-back 
of foliations} more natural. For a morphism of derived DM-stacks
$f : X \longrightarrow Y$, there is an induced $\mathcal{H}$-equivariant
morphism $\mathcal{L}^fX \longrightarrow \mathcal{L}^fY$. For 
a derived foliation $\F \in \Fol(Y)$, realized 
as an $\mathcal{H}$-equivariant derived stack $\F \longrightarrow Y$, 
the pull-back $f^*(\F)$ simply is realized by the pull-back of derived
stacks 
$$f^*(\F)\simeq \F \times_{\mathcal{L}^fY}\mathcal{L}^fX,$$
equipped with its natural projection down to $\mathcal{L}^fX$. 

\subsection{Examples} 

We finish this Section by giving some classes of examples of derived foliations.

\subsubsection{Lie algebroids.}\label{liealgebroid}
Let us now assume that $X$ is a smooth DM stack. Its tangent 
sheaf $\T_X$ is a sheaf of Lie-algebras on the small \'etale site $X_{et}$
where the Lie bracket is the usual bracket of vector fields.
Recall from \cite{nuit} that a Lie algebroid on $X$ consists of a 
vector bundle $T$ on $X$, together two additional structures:

\begin{enumerate}
\item a $\C$-linear Lie bracket $[-,-]$ on $T$.

\item an $\OO_X$-linear morphism $a : T \longrightarrow \T_X$.
\end{enumerate}

These data are required to satisfy the following compatibility relation: for any local sections $s,t$ of $T$, and
any function $f$ on $X$
$$[s,ft]=f[s,t]+a(s)(f)t.$$
We can associate to a Lie algebroid on $X$  a natural derived foliation on $X$ 
as follows. We consider the graded $\OO_X$-cdga $Sym_{\OO_X}(T^*[1])$, 
where $T^*$ is the $\OO_X$-linear dual to $T$. The bracket on $T$
induces a $\C$-linear differential $d : T^* \longrightarrow T^* \wedge_{\OO_X}T^*$, 
which endows $Sym_{\OO_X}(T^*[1])$ with the structure of graded mixed
cdga. This is an object in $\Fol(X)$. The cotangent complex of this derived foliation
is $T^*$ by construction, and this the above derived foliation
is obviously smooth. However, it is rigid only when $a : T \longrightarrow \T_X$ 
identifies $T$ with a subbundle of $\T_X$.

It is easy to show that this construction produces a fully faithful $\s$-functor
$$\mathsf{LieAlgbd}/X \longrightarrow \Fol(X)$$
where $\mathsf{LieAlgbd}/X$ is the category of Lie algebroids over $X$. The essential image
of this $\s$-functor can be shown to coincide with the full $\s$-subcategory consisting of smooth derived foliations over 
$X$ (Definition \ref{d4}).
To be more precise, for any vector bundle $V$ on $X$, 
the classifying space of graded mixed structures on 
the sheaf of graded cgda $Sym_{\OO_X}(V[1])$ turns out to be 
discrete and in bijection with Lie algebroid structures on $V^*$. In particular,
we get that the $\s$-category of smooth derived foliations over $Spec\, \C$ 
is equivalent to the usual category of finite dimensional complex Lie algebras.

There is also a relation between derived foliations and \emph{dg-Lie algebroids} as considered in   
\cite{cg,nuit}. To a dg-Lie algebroid $T$ over an affine variety $X=Spec\, A$, 
we can associate its Chevalley-Eilenberg cochain complex 
$C^*(T):=Sym_{A}(T^*[1])$, considered as a graded mixed cdga using the
Lie bracket as mixed structure. Though this will not be relevant in this paper, 
we think that this construction produces
a \emph{fully faithful} $\s$-functor from the full $\s$-subcategory consisting of dg-Lie algebroids that are perfect over $A$
and of amplitude in $[0,\infty)$, to derived foliations over $X$.

\subsubsection{Initial and final foliations.} The $\s$-category $\Fol(X)$ possesses two important distinguished
objects, namely the initial and final objects. The initial object is called the \emph{punctual derived foliation},
or the \emph{trivial derived foliation},
and is denoted by $\mathbf{0}_X \in \Fol(X)$. The corresponding graded mixed cdga $\DR(\F)$ is $\OO_X$
endowed with the trivial graded mixed structure (concentrated in weight and cohomological degree $0$). 
When $X$ is smooth, $0_X$ can  also be represented by
the $0$ Lie algebroids on $X$ and intuitively corresponds to the foliation on $X$ whose
leaves are the points of $X$.

The final object in $\Fol(X)$ is called the \emph{tautological derived foliation}, or 
\emph{the de Rham derived foliation}. It corresponds to the graded mixed cdga $\DR(X)$ which is the 
derived de Rham algebra of $X$. when $X$ is smooth, it can also be represented by 
the tautological Lie algebroid $\T_X$ itself, and intuitively 
is the foliation for which $X$ is its only leaf.

\subsubsection{Algebraically integrable derived foliations.} Suppose that we are given 
a morphism of derived DM-stacks $f : X \longrightarrow Y$ that is  
locally of finite presentation. The relative de Rham algebra of $X$ over $Y$ defines
a sheaf of graded mixed cdga $\DR(X/Y)$ on the small \'etale site 
of $X$, which is a derived foliation over $X$. Its underlying sheaf of graded
cdga's simply is $Sym_{\OO_X}(\LL_{X/Y}[1])$, where 
$\LL_{X/Y}$ is the relative cotangent complex of $X$ over $Y$. This
is called the \emph{derived foliation
induced by the morphism $f$}. We will use the notation
$$\F_f := \DR(X/Y) \in \Fol(X)$$
for this foliation. Note that $\F_f$
can also be understood as the pull-back
$f^*(\mathbf{0}_Y)$, where $\mathbf{0}_Y$ is the punctual foliation described above.

We set the following definition, and use the expression \emph{d-integrable} to avoid confusions
with the usual notions of integrability of singular underived foliations that we will meet later on
(see \ref{malgr}).

\begin{df}\label{alglocint} Let $X$ be a derived DM stack. 
A derived foliation $F$ on $X$ (locally) 
equivalent to one of the form $\F_f=f^*(\mathbf{0}_Y)$, for a (locally defined) morphism $f:X\to Y$ locally of 
finite presentation between derived DM stacks, will be called \emph{algebraically (locally) d-integrable}.
\end{df}

The reason for this name is that the derived foliation $\F_f$ corresponds intuitively to the foliation on $X$
whose leaves on are the derived fibers of the map $f$. See also Remark \ref{malgr}.  
It is obvious to see that $\F_f$ is quasi-smooth (resp. smooth)
if and only if $f$ is quasi-smooth (resp. smooth). Also, 
$\F_f$ is automatically rigid.

\subsubsection{Pfaffian systems as quasi-smooth and rigid derived foliations.} 
Let $X$ be a smooth algebraic variety. Assume that we are
given 
differential forms $w_i \in \Gamma(X,\Omega_X^1)$, for $i=1,\dots,n$, 
such that the graded ideal $(w_1,\dots,w_k) \subset \Gamma(X,Sym_{\OO_X}(\Omega_X^1[1]))$
is stable by the de Rham differential. We chose differential forms 
$w_{ij}\in \Gamma(X,\Omega^1_X)$ such that 
for all $i$ we have
$$d(w_i)=\sum_j w_{ij}\wedge w_j.$$
We assume furthermore that the $k\times k$ matrix of forms $W=(w_{ij})_{ij}$ satisfy the 
following integrability condition
$$d(W) + W\wedge W = 0.$$
Out of these data 
$w_i$ and $W$ as above, we construct a sheaf of graded mixed cdga's on $X$
by considering $Sym_{\OO_X}(\LL[1])$ where $\LL$ is the two terms perfect complex
$$\LL := \left( \xymatrix{\OO_X^k \ar[r]^-{w_*} & \Omega_X^1} \right).$$ 
The graded mixed structure on $Sym_{\OO_X}(\LL[1])$ is itself 
determined by a morphism of complexes of sheaves of $\C$-vector spaces
$$\LL \longrightarrow \wedge^2_{\OO_X}\LL$$
compatible with the de Rham differential on $\Omega_X^1$. 
Such a morphism is obtained for instance by specifying 
a morphism $\OO_X^k \longrightarrow \OO_X^k \otimes_{\OO_X} \Omega_X^1$
which is a flat connection on the vector bundle $\OO_X^k$. Therefore, 
the matrix $W$ defines such a graded mixed structure, and therefore
a derived foliation on $X$.

The derived foliations defined above depends not only on the $w_i$'s, but also on the
choice of the matrix $W$. It is clear that such derived foliations are quasi-smooth
and rigid. We call such derived foliations
\emph{Pfaffian derived foliations} for obvious reasons. Derived foliations which are 
algebraically d-integrable are always locally, for the Zariski topology, equivalent to 
Pfaffian derived foliations.  

\subsubsection{Pull-backs of smooth and rigid foliations.} 
Let $f : X \longrightarrow Y$ be a morphism of smooth algebraic varieties and 
$\F \in \Fol(Y)$. We have seen that there is a pull-back $f^*(\F) \in \Fol(X)$. 
It is easy to see that when $\F$ is smooth, then $f^*(\F)$ is always quasi-smooth.
If $\F$ is moreover rigid, then so is $f^*(\F)$. 
We will see later that, at least if one admits $Y$ to be a formal scheme, 
all rigid and quasi-smooth derived foliations are locally of this form (see Proposition \ref{p4}).
This follows from an important property of pull-backs, namely the existence of
a homotopy push-out of cotangent complexes
$$\xymatrix{f^*(\Omega_Y^1) \ar[r] \ar[d] & \Omega_X^1 \ar[d] \\
f^*(\LL_\F) \ar[r] & \LL_{f^*(\F).}
}$$
The intuition behind this is that the leaves of $f^*(\F)$ are obtained by 
derived pull-back along $f$ of the leaves of $\F$. The defect of tranversality of $f$ with 
respect to the leaves of $\F$ implies that the leaves of $f^*(\F)$ can be singular, but 
are always quasi-smooth.

\subsubsection{Derived foliations, truncations and singular algebraic foliations.}\label{truncation} 
Our derived foliations bear an important relation with the \emph{singular foliations}
classically considered in the algebraic and holomorphic contexts. There are several possible
definitions of singular foliations in the literature. 
In \cite[Def. 1.11]{ba} they are defined as full coherent differential ideals of the sheaf of differential $1$-forms.
For our purposes, the property of being full will be irrelevant, and will simply define 
a \emph{singular foliation} (on a smooth variety $X$)  as a coherent subsheaf 
$D \subset \Omega_X^1$ satisfying the condition
$$d(D) \subset D\wedge \Omega_X^1 \subset \Omega_X^2,$$
where $d : \Omega_X^1 \longrightarrow \Omega_X^2$ is the 
de Rham differential. This is equivalent to \cite[Def. 6.1.1]{ay}.

Let now $\F \in \Fol(X)$ be a derived foliation on a smooth variety $X$. We consider the
anchor map $a : \Omega_X^1 \longrightarrow \LL_{\F}$ and its induced morphism 
$\Omega_X^1 \longrightarrow H^0(\LL_{\F}).$
We let $D \subset \Omega_X^1$ be the kernel of the above map, which is a coherent subsheaf 
of $\Omega_X^1$. As $a$ comes from a morphism $\DR(X) \longrightarrow Sym_{\OO_X}(\LL_{\F}[1])$ of graded mixed cdga's over $X$,
it follows that the ideal in $\Omega_X^*$ generated by $D$ is in fact 
a differential ideal ($d(D) \subset D\wedge \Omega_X^1$), i.e. $D$ is an underived singular foliation on $X$.
\begin{df}\label{trunc} The kernel $D$ of $H^0(a):\Omega_X^1 \longrightarrow H^0(\LL_{\F})$  is a singular foliation
on $X$ called the \emph{truncation of} $\F$, and denoted by $\tau_{0}(\F)$. 
\end{df}
This produces an $\s$-functor $\tau_{0}$ 
from $\Fol(X)$ to the category of singular foliations on $X$. However, we will see later that 
not all singular foliation arise this way, and the existence of a derived enhancement of a singular foliation is a sublte question related to local integrability. This question will be
studied in details for quasi-smooth and rigid derived foliation (see Corollary \ref{cp6}). 

Conversely, let $D \subset \Omega_X^1$ be a singular foliation on $X$. We can construct 
a graded algebra 
$$\DR(D):=\bigoplus_{p} \left(\Omega_X^p/<D>\right)[p],$$
which is the (underived) quotient of the algebra of differential forms $\Omega_X^*$ by the 
graded ideal generated by $D$. As $D$ is a differential ideal, the de Rham differential
induces a graded mixed structure on $\DR(D)$ in such a way that the canonical morphism
$$\DR(X) \longrightarrow \DR(D)$$
becomes a morphism of graded mixed cdga's. Note that however, $\DR(D)$ does not satisfy the condition
of definition \ref{d1}, as $Sym(\Omega^1_X/D[1])$ involves derived 
wedge powers of the coherent sheaf $\Omega^1_X/D[1]$, and might differ from
$(\Omega^p/<D>)[p]$. Therefore, $\DR(D)$ 
\emph{does not} define a derived foliation in our sense. To be more precise, 
the underlying graded algebra $\DR(D)$ is of the form
$Sym_{\OO_{X}}^u(\Omega^1_X/D[1])$, where $Sym^u$ is the underived symmetric algebra functor. The
construction $D \mapsto \DR(D)$ is easily seen to be an equivalence of categories, from
singular foliations on $X$ to graded mixed algebras of the form $Sym_{\OO_{X}}^u(M[1])$ such that 
$\Omega_X^1 \longrightarrow M$ is a coherent quotient. This however does \emph{not} define 
a functor from singular foliations to derived foliations. 

Finally, the truncation $\s$-functor $\tau_0$ enjoys a certain universal property described as follows. 
The derived foliation $\F$ has a realization 
$|\DR(\F)|$, and in the same way, the truncation $\tau_0(\F)$ has  
an underived de Rham complex
$(\Omega_X^*/<D>)$ where the differential is the de Rham differential. There is an induced canonical
morphism
$$|\DR(\F)| \longrightarrow (\Omega_X^*/<D>).$$
This morphism is far from being an isomorphism in general, although 
it is an isomorphism in low degrees under appropriate conditions (see Proposition \ref{p9}).\\

Note that there is a notion of pullback of singular foliation along an map $f:X\to Y$, with $X$ and $Y$ smooth. 
For simplicity we state it in the algebraic case, but the analytic and formal versions are alike. 
If $D\subset \Omega^1_Y$ is an algebraic singular foliation on $Y$, 
then it's pullback is, by definition, the subsheaf image of the composite map 
$f^*(D) \to f^*(\Omega^1_Y) \to \Omega^1_X$. Using this notion we
give the following definition.

\begin{df}\label{locintsingfol}
An algebraic (resp. formal, resp. analytic) singular foliation on a smooth algebraic variety (resp. formallly 
smooth formal scheme, resp. smooth analytic space) $X$ is \emph{locally integrable} if locally in the Zariski 
topology (resp. locally formally, resp. locally in the analytic topology) at each point of $X$ there exists 
a morphism $f : X \rightarrow Y$ such that $D$ equals is the pull-back by $f$ of the punctual foliation on $Y$.
Equivalently, $D$ is the image of $f^*(\Omega_Y^1) \rightarrow \Omega_X^1$.
\end{df}

Clearly, if a derived foliation $\F$ is d-integrable in the sense of \ref{alglocint}, then 
its truncation $\tau_0(\F)$ is integrable in the sense above.

\begin{rmk}\label{malgr}\emph{In the literature, a (locally) integrable singular foliation  is also sometimes 
called (locally) \emph{completely integrable}. 
We note here that 
the notion of integrability of Definition \ref{locintsingfol} only coincides with the notion of
\cite[\S 3]{mal2} if one further assume that $f$ is generically smooth on $X$ (this is condition $b)$ at p. 73 of 
loc.cit). Therefore, to distinguish between the two, we will refer to 
Malgrange's stronger notion as \emph{local strong integrability.}}
\end{rmk}

\subsection{Foliations over formal completions}
 
Let 
$X$ be a smooth affine variety and $Y \subset X$ be 
a closed subvariety. For simplicity we assume that 
the ideal of $Y$ in $X$ is generated by a regular sequence
$(f_1,\dots,f_k)$.
We denote by 
$\widehat{\OO}_Y$ the ring of functions on the formal completion
of $X$ along $Y$. 

Recall that we have a module $\widehat{\Omega_Y^1}$ of differential forms
on the formal completion of $X$ along $Y$, defined as the formal 
completion of $\Omega_X^1$ along (the ideal defining) $Y$. This is a
locally free $\widehat{\OO}_Y$-module of rank the dimension of $X$. Moreover, it comes
equipped with a canonical derivation $d : \widehat{\OO_Y} \longrightarrow 
\widehat{\Omega_Y^1}$ which extends to a full structure of  
graded mixed cdga $\widehat{\DR}(Y)$ on the graded cdga
$Sym_{\widehat{\OO}_Y}(\widehat{\Omega_Y^1}[1]).$
We can then consider the $\s$-category of 
graded mixed cdga's $A$ endowed with a morphism
$$u : \widehat{\DR}(Y) \longrightarrow A$$
and satisfying the following conditions.

\begin{itemize}
\item The morphism $u$ induces an isomorphism $\widehat{\OO}_Y \longrightarrow A(0)$.

\item The $\widehat{\OO}_Y$-dg-module $A(1)[-1]$ is perfect
and connective.

\item The natural morphism of graded cdga's
$$Sym_{\widehat{\OO}_Y}(A(1)) \longrightarrow A$$
is a quasi-isomorphism.

\end{itemize}

Let us denote by $\widehat{\Fol}(\widehat{Y})$ 
the opposite $\s$-category of the above $\s$-category of graded
mixed cdga's \emph{under} $\widehat{\DR}(Y)$. We call its objects
\emph{formal derived foliations on the formal scheme $\widehat{Y}$}. 
On the other hand, we may identify the formal scheme $\widehat{Y}$ with
its associated (derived) stack, and therefore we may consider its $\s$-category
of foliations $\Fol(\widehat{Y})$, according to Definition \ref{d2}.

\begin{prop}\label{p3}
There exists a natural equivalence of $\s$-categories
$$\widehat{\Fol}(\widehat{Y})\simeq \Fol(\widehat{Y}).$$
\end{prop}

\noindent\textbf{Proof.} 
The formal completion of $X$ along $Y$, denoted by $\widehat{Y}$, 
is equivalent, as an object of $\dSt$, to a colimit 
$$\widehat{Y}:=colim_n Y_n$$
where $Y_n \subset X$ is the closed sub-scheme defined by the
equations $(f_1^n,\dots,f_k^n)$, where $(f_1,\dots,f_k)$ is the regular sequence
generating the ideal of $Y$ in $X$. We thus have
$\Fol(\widehat{Y})\simeq lim_n \Fol(Y_n)$.
Now, the right hand side is directly related to $\widehat{\Fol}(\widehat{Y})$ by the limit $\s$-functor
$$lim : lim_n \Fol(Y_n) \longrightarrow \widehat{\Fol}(\widehat{Y}).$$
Note that this is well defined as $lim_n\DR(Y_n)\simeq \widehat{\DR}(Y)$, 
because $\Omega^1_{\widehat{Y}} \simeq lim_n \mathbb{L}_{Y_n}$.
The inverse functor is defined by sending 
a graded mixed cgda $A$ under $\widehat{\DR}(Y)$ to 
its families of restriction 
$$\{A\otimes_{\widehat{\DR}(Y)}\DR(Y_n)\}_n \in lim_n \Fol(Y_n).$$
The fact that these two $\s$-functors are inverse to each others
is an immediate consequence of the fact that the natural $\s$-functor
$$\Perf(\widehat{Y}) \longrightarrow lim_n \Perf(Y_n)$$
is an equivalence.
\hfill $\Box$ \\

The following is the formal version of Definition \ref{alglocint}.

\begin{df}\label{forlocint}
Let $\widehat{Y}$ (respectively, $\widehat{Y')}$) be the formal completion of a smooth affine scheme $Y$ (resp. $Y'$) along an ideal generated by a regular sequence. We say that a foliation $\F$ on $\widehat{Y}$ is \emph{formally (locally) d-integrable} if there exists a (locally defined) morphism of formal schemes $f:\widehat{Y} \to \widehat{Y'}$ such that $\F$ is (formally locally) equivalent to the pullback $f^*(\mathbf{0}_{\widehat{Y'}})$  of the final foliation $\mathbf{0}_{\widehat{Y'}}$ on $\widehat{Y'}$ via $f$.

If $X$ is a smooth variety, $x\in X$, and $\widehat{X}_x$ denotes the formal completion of $X$ at $x$, then an algebraic foliation $\F \in \F ol(X)$ is said to be \emph{formally d-integrable at} $x$ if its restriction $\widehat{\F}$ to $\widehat{X}_x$ (which is a formal foliation on $\widehat{X}_x$) 
is formally d-integrable according to the previous definition.
\end{df}

Note that, in the above definition, the underlying graded
cdga's of  $f^*(\mathbf{0}_{\widehat{Y'}})$ is $Sym_{\OO_X}(\widehat{\LL}_{\widehat{Y}/\widehat{Y'}}[1])$, where 
$\widehat{\LL}_{\widehat{Y}/\widehat{Y'}}$ is the cotangent complex of $f:\widehat{Y} \to \widehat{Y'}$.

\subsection{Formal structure of rigid quasi-smooth derived foliations}

In the proposition below we let $\widehat{\mathbb{A}}^n$ be the formal
completion of $\mathbb{A}^n$ at $0$.

\begin{prop}\label{p4}
Let $X$ be smooth variety and $\F \in \Fol(X)$ be a rigid
and quasi-smooth derived foliation on $X$. Then, Zariski locally 
on $X$, there exists $n\geq 0$ and a smooth and rigid derived foliation 
$\F'$ on the formal scheme $X\times \widehat{\mathbb{A}}^n$, such that 
$\F$ is the pull-back of $\F'$ by the zero section 
$X \longrightarrow X \times \widehat{\mathbb{A}}^n$.
\end{prop}

\noindent\textbf{Proof.} 
We will freely use the material and notations about internal De Rham complexes and their realizations, recalled in Section \ref{reminder}. First of all, the statement being Zariski local, we may 
assume that $X=Spec\, A$ is a smooth affine variety. 
We consider $\F \in \Fol(X)$, a rigid and 
quasi-smooth derived foliation, which corresponds to a 
graded mixed cdga $\DR(\F)$ with an identification $\DR(\F)(0)\simeq A$ and
satisfying the conditions of definitions \ref{d1} and \ref{d4}. There is a natural morphism
$\Omega_A^1 \longrightarrow \mathbb{L}_\F,$
whose cone, by the rigidity and quasi-smooth hypothesis, must be 
of the form $N^*[1]$, for a vector bundle $N^*$
over $X$. By localizing further on $X$ we may
suppose $N^*$ isomorphic to $\OO_X^n$ for some $n\geq 0$.
We consider the push-out of graded mixed cdga's
$$\DR(\F) \longrightarrow \DR(\F)\otimes_{\DR(X)}A.$$
On the underlying graded cdga's, this morphism  looks like
$$Sym_A(\mathbb{L}_{\F}[1]) \longrightarrow Sym_A(N^*[2])$$
and is induced by applying the $Sym$ construction to the boundary map $\LL_\F \rightarrow N^*[1]$. 
Moreover, as $N^*$ is a vector bundle (and $X$ is affine) we see that the graded
mixed structure on the graded mixed cdga $Sym_A(N^*[2])$ must be trivial. 

We now consider the induced morphism on the internal relative 
de Rham algebra
$$\DR^{int}(Sym_A(N^*[2])/\DR(\F)) \longrightarrow \DR^{int}(A/\DR(\F)).$$
We can then consider their internal realization, and the induced
morphism on the internal de Rham cohomology
$$|\DR^{int}(Sym_A(N^*[2])/\DR(\F))| \longrightarrow 
|\DR^{int}(A/\DR(\F))|.$$
This is a new morphism of graded mixed cdga's and thanks to lemma \ref{l1}, 
the right hand side is canonically 
equivalent to $\DR(\F)$, so we get a morphism
$$|\DR^{int}(Sym_A(N^*[2])/\DR(\F))| \longrightarrow \DR(\F).$$
The left hand side is a graded mixed cdga, say $B$, 
whose degree zero part is $|Sym_A(N^*[2])|\simeq \widehat{Sym}_A(N^*)\simeq 
\OO(X\times \widehat{\mathbb{A}}^n)$.
Moreover, by construction, it is not hard to see that 
$B$ is free, as a graded cdga, over the $B(0)$-module $\Omega_A^1\otimes_A B(0)$. Therefore,
Proposition \ref{p3} tells us that
$B$ corresponds to a smooth foliation $\F'$ on $X\times \widehat{\mathbb{A}}^n$.

The morphism of graded mixed cdga's $B \longrightarrow \DR(\F)$ induces a morphism 
$e^*(\F') \longrightarrow \F$ in $\Fol(X)$, where 
$e : X \longrightarrow X \times \widehat{\mathbb{A}}^n$
is the zero section. This morphism is furthermore an equivalence by observing 
the induced morphism on cotangent complexes. In particular, 
$e^*(\F')$ is rigid, and this automatically implies that $\F'$ is rigid as
a foliation on $X\times \widehat{\mathbb{A}}^n$.
\hfill $\Box$ \\

\begin{rmk}\label{r1bis}
\emph{Proposition \ref{p4} has the following conceptual refinement.
First of all, the same proof shows that there exists a globally defined 
pair $(\mathfrak{X}' , \mathcal{F}')$
consisting of a formally smooth formal scheme $\mathfrak{X}'$ 
and a smooth rigid foliation $\mathcal{F}'$ on $\mathfrak{X}'$, together with a map $f:X \to \mathfrak{X}'$
such that $f_{red}$ is an isomorphism and there exists an equivalence $\mathcal{F} \simeq f^*(\mathcal{F}')$.
This is achieved by defining 
$$\mathfrak{X}':= \mathrm{Spf}(|\DR(\F)\otimes_{\DR(X)}A|),$$
and $\mathcal{F}'$ as the derived foliation on $\mathfrak{X}'$ defined by  the graded mixed cdga on $X$
$$|\DR^{int}((\DR(\F)\otimes_{\DR(X)}A)/\DR(\F))|.$$
Note that if $\mathcal{N}^*$ is the vector bundle on $X$ such that $\mathcal{N}^*[1]$ is the cofiber of the canonical map
$\Omega_{X}^1 \longrightarrow \mathbb{L}_\F$, the underlying graded cdga $\DR(\F)\otimes_{\DR(X)}A$ 
is $Sym_{\OO_X}(N^*[2])$. However, when $X$ is not affine the mixed structure on 
$Sym_{\OO_X}(N^*[2])$ can possibly be non-trivial. We note here that this mixed structure
captures the difference between the formal scheme $\mathfrak{X}'$ and 
the formal completion $\widehat{\mathcal{N}}$ of the normal bundle of $\mathfrak{X}$ in $\mathfrak{X}'$.
Finally, Proposition \ref{p4} boils down to the fact that, Zariski locally on $X$, $\mathfrak{X}'$ can be chosen of the form
$X\times \widehat{\mathbb{A}}^n$ (and this is simply due that $\mathcal{N}^*$ is locally free of finite rank), and that, under this identification, $f$
becomes identified with the zero section map $X \to X\times \widehat{\mathbb{A}}^n$.
}
\end{rmk}

\begin{rmk}\label{r1}
\emph{Without the rigidity assumption, but still for quasi-smooth 
derived foliations, Proposition \ref{p4} remains true if one
replaces the formal scheme $X\times \widehat{\mathbb{A}}^n$ by 
a formally smooth formal Artin stack containing $X$ as its reduced part. }
\end{rmk}

When restricting to the formal completion at a point $x \in X$, Proposition \ref{p4} has
the following important consequence. 

\begin{cor}\label{cp4}
Let $\F \in  \Fol(X)$ be a derived quasi-smooth and rigid derived foliation over a smooth
algebraic variety $X$. Let $x \in X$ and
$\widehat{X}_x=Spf(\widehat{\OO}_{X,x})$ be the formal completion of $X$ at $x$, and $\widehat{\F} \in \Fol(\widehat{X}_x)$
the restriction of $\F$. Then, there
exists $m\geq 0$, a morphism $f : \widehat{X}_x \longrightarrow \widehat{\mathbb{A}}^m$ and an isomorphism 
$$f^*(\mathbf{0}_{\widehat{\mathbb{A}}^m}) \simeq \widehat{\F}.$$
In other words, $\F$ is formally d-integrable at each point (Def. \ref{forlocint}). 
\end{cor}

\noindent\textbf{Proof.} Indeed, by Proposition \ref{p4}, $\F$ is, locally at $x$, the pull-back of 
a smooth and rigid derived foliation $\F'$ on $X \times \widehat{\mathbb{A}}^n$. Thus, $\widehat{\F}$ becomes
isomorphic to the pull-back of a smooth and rigid derived foliation $\widehat{\F}'$ on
$\widehat{X}_x \times \widehat{\mathbb{A}}^n$. By the formal version of Frobenius theorem (see for instance
\cite[Thm. 2]{bol}), 
we know that $\widehat{\F}'$ is formally integrable and thus d-integrable as it is a smooth
derived foliation. This implies that $\widehat{\F}$ is d-integrable.
\hfill $\Box$ \\

\begin{rmk}
\emph{The above corollary is also true for non-quasi-smooth derived foliation, but 
$\widehat{\mathbb{A}}^m$ must be replaced by a more general, eventually not formally smooth, 
derived formal schemes. It shows in particular that not all singular foliation is the truncation
of a quasi-smooth and rigid derived foliation, not even at the formal level. Indeed, 
formal integrability is not always satisfied for singular foliations}.
\end{rmk}

\subsection{Leaves of a derived foliation}

Given a derived foliation $\F$ on a smooth variety $X$, and a point $x\in X$, it is possible to define
the notion of \emph{leaf of $\F$ passing through $x$}, at least at the formal level. For this, we start
by recalling the notion of formal moduli problems in the sense of \cite{lupmf} (see also 
\cite{tobour} for an overview). 

We let $\dgart^*$ be the $\s$-category of commutative, Artinian, connective, local and augmented dg-algebras. 
Recall that these are cdga's $A$, together with an augmentation $A \longrightarrow \C$, and such that 
the following conditions are satisfied:

\begin{itemize}
\item The ring $H^0(A)$ is local and artinian.

\item The vector space $H^*(A)$ is finite dimensional over $\C$.

\end{itemize}

By definition, a formal moduli problem is an $\s$-functors
$$F : \dgart^* \longrightarrow \mathbf{Top}$$
satisfying the following two conditions:

\begin{itemize}
\item The space $F(\C)$ is contractible.

\item For any cartesian square in $\dgart^*$
$$\xymatrix{A' \ar[d] \ar[r] & B' \ar[d] \\
A \ar[r] & B}$$
such that the map $H^0(A) \longrightarrow H^0(B)$
is surjective, the induced commutative diagram of spaces
$$\xymatrix{F(A') \ar[d] \ar[r] & F(B') \ar[d] \\
F(A) \ar[r] & F(B)}$$
is cartesian in $\mathbf{Top}$.

\end{itemize}

Graded mixed cdga's can be used to define formal moduli problems as follows
(see also \cite{cptvv} for more about the relations between formal geometry and
graded mixed cdga's). To start with, there
is an $\s$-functor
$$\DR(\C/-) : \dgart^* \longrightarrow \medg$$
sending an augmented cdga $A \rightarrow \C$ to $\DR(\C/A)$, the de Rham algebra of $\C$ relative to $A$, through the augmentation.
This $\s$-functor is easily seen to be fully faithful, the artinian cdga $A$ being recovered
as the realization $|\DR(\C/A)|$.
According to the classification of formal moduli problems by dg-Lie algebras, this $\s$-functor
can also be described as follows. An object $A \in \dgart^*$ gives rise to a tangent 
dg-lie algebra $\ell_A:=\T_{\C/A}$, the dg-lie algebra of derived $A$-linear derivations on $\C$. The
graded mixed cdga $\DR(\C/A)$ is then canonically equivalent to $C^*(\ell_A)$, the Chevalley complex
of the dg-Lie algebra $\ell_A$, considered as a graded mixed cdga (see \cite{tobour} for details). 

For a given graded mixed cdga $B$ we thus define its formal spectrum
$\Spf\, B : \dgart^* \longrightarrow \mathbf{Top}$
by 
$$(\Spf\, B)(A):=\mathrm{Map}_{\medg}(B,\DR(\C/A)).$$
The $\s$-functor $\Spf\, B$ is then a formal moduli problem, called the \emph{formal spectrum} of $B$.\\

Let now $\F$ be a derived foliation on the smooth variety $X$, and $x\in X$. Taking an affine chart around $x$ we can
consider $X$ to be affine, say $X=Spec\, A$.
We start by considering the augmentation induced by  the point $x$
$$\DR(\F) \longrightarrow A \longrightarrow \C.$$
This is a morphism of graded mixed cdga's and we can form the corresponding 
push-out
$\DR(\F_x):=A \otimes_{\DR(\F)}\C.$
This is a new graded mixed cdga which corresponds to a derived foliation over 
$Spec\, \C$, and whose cotangent complex is $\LL_{\F,x}[1]$, the shift of the
fiber at $x$ of the cotangent complex of $\F$. The formal spectrum of $\DR(\F_x)$ 
defines a formal moduli problem
$$\widehat{\F}_x := \Spf\, (\DR(\F_x)) : \dgart^* \longrightarrow \mathbf{Top}.$$

\begin{df}\label{dleaf}
The formal moduli problem $\widehat{\F}_x$ defined above is called
the \emph{formal leaf of $\F$ passing through $x$}.
\end{df}

Note that the tangent dg-Lie algebra of the formal leaf $\widehat{\F}_x$ 
is given by $\T_{\F,x}[-1]$, the fiber at $x$ of the tangent complex of 
the derived foliation $\F$. In particular, amplitude considerations tell us that 
$\widehat{\F}_x$ is always representable by a derived formal scheme (see \cite{lupmf}). 

When applied to the tautological foliation, $\DR(\F)=\DR(X)$, 
we get that $\DR(\F_x)=Sym_{\C}(\Omega^1_{X,x}[1])$ with the trivial
mixed structure. In this case we see that $\Spf\, (\DR(\F_x))$ is canonically
isomorphic to $\widehat{X}_x$, the formal completion of $X$ at $x$.
Moreover, using the canonical map $\DR(X) \longrightarrow
\DR(\F)$, the formal moduli $\widehat{\F}_x$ comes equipped with a canonical
morphism 
$$\widehat{\F}_x \longrightarrow \widehat{X}_x$$
to the formal completion of $X$ at the point $x$. When $\F$ is rigid this morphism
is a closed embedding of derived formal schemes, and thus $\widehat{\F}_x$ sits inside
$\widehat{X}_x$ as a closed formal derived subscheme. 

The morphism induced by $\widehat{\F}_x \longrightarrow \widehat{X}_x$ at the level of rings of functions, can 
be described as follows. The commutative dg-algebra of functions on $\widehat{F}_x$ is
$|\DR(\F_x)|$, the realization of $\DR(\F_x)$. The natural morphism
$\DR(X) \longrightarrow \DR(\F_x)$ induces a canonical morphism of cdga's
$$\OO(\widehat{X}_x)\simeq |Sym_{\C}(\Omega^1_{X,x}[1])| \longrightarrow
|\DR(\F_x)|=:\OO(\widehat{\F}_x),$$
which is the induced morphism obtained by restriction of functions along the embedding 
$\widehat{\F}_x \subset \widehat{X}_x$.\\

Note that one can define \emph{analytic} and \emph{algebraic} leaves as well, but of course, contrary to the formal case, their
existence is not guaranteed. We give below the definition in the algebraic setting, the case of 
analytic leaves being completely similar. 

We define an \emph{algebraic leaf} of $\F$ on $X$ to be a derived scheme $L$ together with 
a pointwise injective morphism of derived schemes
$j : L \longrightarrow X$
such that the following condition is satisfied: for all point $x \in L$, 
there exists an equivalence of formal derived schemes $\widehat{L}_x \simeq \widehat{F}_x$, 
that makes the diagram below commutative
$$\xymatrix{
\widehat{L}_x \ar[dr]_-{j} \ar[rr]^\sim & & \widehat{F}_x \ar[dl] \\
& \widehat{X}_x. &
}$$

When $\F$ is globally integrable by a morphism $f : X \longrightarrow Y$, then clearly 
the inclusions of derived fibers $f^{-1}(y) \longrightarrow X$ are 
leaves in the sense above. As we will see later (Corollary \ref{cp6'}), under an appropriate codimension $2$ condition,
rigid and quasi-smooth derived foliations always admit \emph{analytic leaves} locally. i.e. their formal 
leaves are in fact "convergent".

\section{The analytic theory}

The general notion of derived foliations has  
a complex analytic analogue. We will not need the most  
general definition, that would require some advanced tools of derived analytic geometry (see e.g. \cite{porta}), 
and we will restrict ourselves
to derived foliations over smooth complex analytic spaces, for which the basic definitions
can be given more directly.

\subsection{Analytic derived foliations}

Let $X$ be a smooth complex analytic space. It has a
sheaf of holomorphic $1$-forms $\Omega_X^1$, and 
a de Rham algebra $\DR(X):=Sym_{\OO_X}(\Omega_X^1[1])$. This
is a sheaf of graded cdga's over $X$, which is equipped with a canonical
graded mixed structure given by the holomorphic de Rham differential. 

\begin{df}\label{d5}
A \emph{holomorphic or analytic derived foliation over $X$} consists
of a sheaf $A$ of graded mixed cdga's over $X$, together with
a morphism of sheaves of graded mixed cdga's
$$\DR(X) \longrightarrow A$$ satisfying the following conditions:
\begin{enumerate}
\item The induced morphism $\OO_X \longrightarrow A(0)$ is 
a quasi-isomorphism.

\item The complex of $\OO_X$-modules $A(1)[-1]$ is perfect
and connective.

\item The natural morphism $Sym_{\OO_X}(A(1)) \longrightarrow A$
is a quasi-isomorphism of sheaves of graded cdga's.

\end{enumerate}
\end{df}

The analytic derived foliations over a complex manifold $X$ 
form an $\s$-category, denoted by $\Fol(X)$. It is a full sub-$\s$-category
of the $\s$-category of sheaves of graded mixed $\DR(X)$-algebras over $X$. 
For any morphism $f : X \longrightarrow Y$ of complex manifolds, we have a pull-back $\s$-functor
$$f^* : \Fol(Y) \longrightarrow \Fol(X).$$
It is defined as in the algebraic case. There is a natural morphism $f^{-1}(\DR(Y)) \longrightarrow \DR(X)$ of sheaves
of graded mixed cgda's on $Y$. For $\F \in \Fol(Y)$, 
corresponding to a sheaf of graded mixed cdga $\DR(\F)$,  we define 
$f^*(\F) \in \Fol(X)$ as the derived foliation associated to the
sheaf of graded mixed cdga's given by the base change
$$\DR(f^*(\F)):=\DR(X) \otimes_{f^{-1}(\DR(Y))}f^{-1}(\DR(\F)).$$

As in Definition \ref{d4}, we have the notions of \emph{smooth, quasi-smooth, and rigid
derived foliations} over a complex manifold. As in Definition \ref{trunc}, we have a notion of
\emph{truncation} of an analytic derived foliation on a complex manifold; this truncation is an analytic singular foliation on the same complex space.\\

The following is the analytic version of Definition \ref{alglocint}.

\begin{df}\label{analocint}
An analytic derived foliation $\F$ on a complex manifold $X$ is \emph{(locally) d-integrable} if there exists a 
(locally defined) analytic map $F:X\to Y$ of complex manifolds and a (local in the analytic topology) equivalence 
$\F \simeq f^*(\mathbf{0}_{Y})$, where $\mathbf{0}_Y$ is the final derived analytic foliation on $Y$.  
\end{df}

\subsection{Analytification}\label{analitificazione}

Let $X$ be a smooth algebraic variety and $X^h$ be the corresponding
complex analytic space. We are going to construct an analytification $\s$-functor
$$(-)^h : \Fol(X) \longrightarrow \Fol(X^h).$$
We have a morphism of ringed spaces
$$u : (X^h,\OO_X^h) \longrightarrow (X,\OO_X).$$
This morphism induces a canonical isomorphism $u^*(\Omega_X^1) \simeq \Omega_{X^h}^1$ 
of vector bundles on $X^h$. This extends to a natural isomorphisms
$u^*(\Omega_X^*)\simeq \Omega_{X^h}^*$, which is compatible with the de Rham 
differential in the sense that the composed morphism
$$u^{-1}(\DR(X)) \longrightarrow u^*(\DR(X)) \simeq \DR(X^h)$$
is not only a morphism of graded cdga's but of graded mixed cdga's.

For an algebraic derived foliation $\F \in \Fol(X)$, corresponding to 
a morphism of sheaf of graded mixed cdga's $\DR(X) \longrightarrow \DR(\F)$, 
we consider
$$\DR(X^h) \longrightarrow u^{-1}(\DR(\F)) \otimes_{\DR(X)}\DR(X^h).$$
This defines a derived foliation $\F^h \in \Fol(X^h)$. Obviously, the
construction $\F \mapsto \F^h$ is functorial and defines
an $\s$-functor 
$$(-)^h : \Fol(X) \longrightarrow \Fol(X^h).$$

\begin{df}\label{d6}
The \emph{analytification $\s$-functor for derived foliations} is the
$\s$-functor 
$$(-)^h : \Fol(X) \longrightarrow \Fol(X^h)$$
defined above.
\end{df}

The analytification $\s$-functor shares the following straightforward 
properties. 

\begin{itemize}
\item (Functoriality) Let $f : X \longrightarrow Y$ be a morphism 
of smooth algebraic varieties and $f^h : X^h \longrightarrow Y^h$
the corresponding morphism of complex spaces. Then, we have a naturally 
commutative diagram of $\s$-functors
$$\xymatrix{
\Fol(Y) \ar[r]^-{f^*} \ar[d]_-{(-)^h} & \Fol(X) \ar[d]^-{(-)^h} \\
\Fol(Y^h) \ar[r]_-{(f^h)^*} & \Fol(X^h).}$$

\item A derived foliation $\F \in \Fol(X)$ is smooth (resp. quasi-smooth, resp.
rigid) if and only if $\F^h$ is smooth (resp. quasi-smooth, resp. rigid).

\item For any smooth variety $X$, the
analytification $\s$-functor $(-)^h : \Fol(X) \longrightarrow \Fol(X^h)$
is conservative.

\item Exactly as done in \ref{truncation} for the algebraic case, there is a truncation $\s$-functor $\tau_0: \Fol (Y) \to \mathsf{SingFol}(Y)$ from analytic derived foliations to analytic singular foliations over a complex manifold $Y$. Moreover, it is easy to check that, if $X$ is a smooth algebraic variety, the following diagram commutes $$\xymatrix{\Fol (X) \ar[d]_-{\tau_0} \ar[r]^-{(-)^h} & \Fol(X^h) \ar[d]^-{\tau_0} \\ \mathsf{SingFol}(X) \ar[r]_-{(-)^h} & \mathsf{SingFol}(X^h) .}$$

\end{itemize}

When $X$ is smooth and proper GAGA implies furthermore the following
statement.

\begin{prop}\label{p5}
For a smooth and proper algebraic variety $X$ the analytification $\s$-functor
$$(-)^h : \Fol(X) \longrightarrow \Fol(X^h)$$
is an equivalence. 
\end{prop}

\noindent\textbf{Proof.} We let $\Omega_X^*[*]:=Sym_{\OO_X}(\Omega_X^1[1])$ be the sheaf of cdga's on $X$.
Its analytification is $\Omega_{X^h}^*[*]=Sym_{\OO_{X^h}}(\Omega_{X^h}^1[1])$. The analytification
functor produces an dg-functor between dg-categories of dg-modules
$$\Omega_X^*[*]-Mod \longrightarrow \Omega_{X^h}^*[*]-Mod.$$
We restrict this $\s$-functor to the full sub-dg-categories
of perfect dg-modules, i.e. sheaves of dg-modules which locally 
are obtained by finite limits and colimits of $\Omega_X^*[*]$ (resp. 
$\Omega_{X^h}^*[*]$), and pass to ind-completion to get 
$$\mathsf{IndPerf}(\Omega_X^*[*]) \longrightarrow \mathsf{IndPerf}(\Omega_{X^h}^*[*]).$$
By GAGA this is an equivalence.
We recall here that for any graded mixed cdga $A$, the dg-category 
of dg-modules $A-Mod\simeq \mathsf{IndPerf}(A)$ has a canonical action of the
group stack $\mathcal{H}=B\mathbb{G}_a \rtimes \mathbb{G}_m$. 
The $\mathbb{G}_m$-action is induced by the grading on $A$, while the
action of $B\mathbb{G}_a$ by the mixed structure. The dg-category of 
fixed points by $\mathcal{H}$ is moreover equivalent to the 
dg-category of graded mixed $A$-dg-modules (see \cite{pato}).
By sheafification, this
implies that the group $\mathcal{H}$ acts on both dg-categories
$\mathsf{IndPerf}(\Omega_X^*[*])$ and $\mathsf{IndPerf}(\Omega_{X^h}^*[*])$ and the analytification dg-functor
becomes an $\mathcal{H}$-equivariant dg-equivalence
$$\mathsf{IndPerf}(\Omega_X^*[*]) \simeq \mathsf{IndPerf}(\Omega_{X^h}^*[*]).$$
We apply the fixed points construction (see \cite{pato} for details) and get this way a
new equivalence of $\s$-categories
$$\mathsf{IndPerf}(\Omega_X^*[*])^{h\mathcal{H}} \simeq \mathsf{IndPerf}(\Omega_{X^h}^*[*])^{h\mathcal{H}}.$$
The analytification functor being compatible with tensor products, 
the above $\s$-functor has a natural symmetric monoidal structure and thus induces an $\s$-equivalence
on the level of $\s$-categories of commutative algebras. The proposition follows
by the observation that $\Fol(X)$ (resp. $\Fol(X^h)$) is a full sub-$\s$-category
of the $\s$-category of commutative algebras in $\mathsf{IndPerf}(\Omega_X^*[*])^{h\mathcal{H}}$
(resp. in $\mathsf{IndPerf}(\Omega_{X^h}^*[*])^{h\mathcal{H}}$) and that these sub-$\s$-categories
match by the above equivalence.
\hfill $\Box$ \\

\subsection{Analytic integrability}

We have seen that quasi-smooth and rigid derived foliations
are always formally d-integrable, a property which distinguishes them from
the underived singular foliations. We now study analytic d-integrability (Def. \ref{analocint}) of quasi-smooth and
rigid derived foliations, locally in the analytic 
topology. 
We think it is not true that analytic d-integrability always holds for quasi-smooth and rigid derived foliations, 
but we will see below (Proposition \ref{p6} and Corollary \ref{cp6'}) that they are always 
locally integrable
under a rather common codimension $\geq 2$ condition. \\

Let $\F \in \Fol(X)$ be a quasi-smooth and rigid derived foliation
on a smooth algebraic variety $X$ and $\F^h \in \Fol(X^h)$ its
analytification. We consider the truncation $\tau_0(\F)$ (Definition \ref{trunc}), which is an algebraic singular foliation on $X$, and its
analytification $\tau_0(\F^h)$, which is an analytic singular foliation
on $X^h$. 

The cotangent complex $\LL_{\F^h}$ is perfect complex of amplitude $[-1,0]$ on $X^h$.

\begin{df}\label{dsing}
With the above notations and assumptions on $\F$, the \emph{smooth locus of $\F$}
is the Zariski open subset in $X$ of points where $\LL_F$ is quasi-isomorphic to a vector bundle sitting in degree $0$.
Its closed complement $Sing(\F) \subset X$ is called the \emph{singular locus of $\F$}.
\end{df}

Equivalently, since $\F$ is supposed to be quasi-smooth and rigid, 
$Sing(\F)$ is the support of the coherent sheaf $H^1(\T_\F)$, where $\T_\F:=\LL_{\F}^{\vee}$ denotes the tangent complex of $\F$. Note that, in particular, the smooth locus of such an $\F$ might be empty. \\

The following result entails local analytic integrability of the truncation of 
any quasi-smooth and rigid derived foliation as soon as we impose
smoothness outside a codimension $\geq 2$ subset. More precisely, we have the following result.

\begin{prop}\label{p6}
Let $X$ be a smooth irreducible algebraic variety and $\F \in \Fol(X)$ be a quasi-smooth rigid 
derived foliation. We assume that the singular locus $Sing(\F) \subset X$ is at least of codimension
$2$. Then, the truncated analytic singular foliation $\tau_0(\F^h)$ on $X^h$ 
is locally strongly integrable in the analytic topology.
\end{prop}

\noindent\textbf{Proof.} This is in fact an easy consequence of our Corollary \ref{cp4}, which ensures that 
the foliation $\F$ in the statement of Proposition \ref{p6} is formally d-integrable at each point. This implies that its truncation $\tau_0(\F^h)$
is a singular foliation on $X^h$ which is formally integrable at each point. At this point, we would like to
apply \cite[Thm. 3.1]{mal2} that proves that a formally strongly integrable singular foliation is analytically strongly integrable if its singular locus has codimension $\geq 2$ (for the notion of formally or analytically strongly integrable singular foliation, see Remark \ref{malgr}) to deduce that $\tau_0(\F^h)$ is, a fortiori, analytically integrable locally around
each point of $X$. But in order to do this, we need to show that $\tau_0(\F^h)$ is not only formally integrable but also formally strongly integrable (at each point). Now, by hypothesis, the smooth locus of $\F$ is a non-empty Zariski open in $X$, hence dense, since $X$ is irreducible, so our $\F$ is a quasi-smooth rigid derived foliation which is actually smooth on an open Zariski dense subset of $X$. Now, for a quasi-smooth (rigid) derived foliation $\F$ on $X$ which is smooth on a dense open Zariski subset of $X$, formal $d$-integrability of $\F$ at $x\in X$ implies formal strong integrability at $x$ for its truncation $\tau_0(\F)$ (since, for $f: X\to Y$ locally defined at $x$, the fact that the pullback derived foliation $f^*(\mathbf{0}_Y)$ is generically smooth entails generic smoothness for the map $f$ itself\footnote{The reader may easily verify that this implication is false for underived singular foliations.}).  Thus we are in a position to apply \cite[Thm. 3.1]{mal2}, and deduce local analytic (strong) integrability around any $x$ for the analytification $\tau_0(\F^h)$. \hfill $\Box$ \\

An important consequence of Proposition \ref{p6} is the following statement, establishing  a precise relation between
underived singular foliations and quasi-smooth rigid derived foliations.

\begin{cor}\label{cp6}
Let $X$ be a smooth irreducible algebraic variety and $D$ an underived singular foliation on $X$ whose
singular locus $Sing(\F) \subset X$ is of codimension at least 2. Then $D$ is
locally, for the analytic topology, the truncation of a quasi-smooth and rigid derived foliation
if and only if it is formally strongly integrable at each point.
\end{cor}

To finish this section we mention the following stronger version of Proposition \ref{p6}, though it will not be 
used in the rest of the paper.

\begin{cor}\label{cp6'}
Let $X$ be a smooth irreducible algebraic variety and $\F \in \Fol(X)$ be a quasi-smooth rigid 
derived foliation. We assume that the singular locus $Sing(\F) \subset X$ is at least of codimension
$2$. Then $\F^h$ is, locally on $X^h$, a d-integrable derived foliation.
\end{cor}

\textit{Proof.} By Proposition \ref{p6}, it is enough to 
show that, under the codimension $2$ hypothesis in the statement, 
an analytic rigid quasi-smooth derived foliation is d-integrable if its
truncation is strongly integrable. 

We thus let $Y$ be a complex smooth manifold and $\G$ a rigid and quasi-smooth 
derived analytic foliation on $Y$. We assume that $Sing(\G)$, the locus in $Y$ 
where $\mathbb{\LL}_\G$ is not a vector bundle (sitting in degree $0$), is of codimension at least $2$. 
Suppose that the truncation $\tau_0(\G)$ is strongly integrable  by a morphism
$f : Y \longrightarrow S$ (see \ref{malgr}). The morphism $f$ is therefore
automatically generically smooth on $Y$. We fix an isomorphism between 
$\tau_0(\G)$ and the singular foliation induced by $f$.

If $\DR(Y) \longrightarrow \DR(\G)$ is the natural morphism of graded mixed
cdga's over $Y$, we first want to construct a factorization 
$$\xymatrix{\DR(Y) \ar[r] & \DR(Y/S)=f^*(\mathbf{0}_S) \ar[r]^-{\bar{v}} & \DR(\G).}$$
We remind that we have a push-out square of graded mixed cdga's on $Y$
$$\xymatrix{
\DR(Y) \ar[r] & \DR(Y/S) \\
f^{-1}(\DR(S)) \ar[r] \ar[u] & f^{-1}(\OO_S). \ar[u]}$$
This implies that a choice of the factorization as above is equivalent to a choice
of a factorization of the natural morphism $f^{-1}(\DR(S)) \longrightarrow  \DR(\G)$ as
$$\xymatrix{f^{-1}(\DR(S)) \ar[r] & f^{-1}(\OO_S) \ar[r]^-{v} & \DR(\G).}$$
By the universal property of the de Rham graded mixed cdga $\DR(S)$, choosing such a factorization
is equivalent to promote the natural morphism of graded cdga's $f^{-1}(\OO_S) \longrightarrow 
\DR(\G)$
to a morphism of graded mixed cdga's, where the mixed structure on the left hand side is trivial.

We consider the sheaf of algebras $H^0(|\DR(\G)|)$. Its comes equipped with a canonical morphism
of sheaves of cdga's $H^0(|\DR(\G)|) \longrightarrow |\DR(\G)|$. Moreover, by adjunction, there is a
natural morphism of graded mixed cdga's
$$|\DR(\G)| \longrightarrow \DR(\G)$$
where the graded mixed structure on the left hand side is the trivial one. We thus obtain this way 
a morphism of graded mixed cgda's
$$u : H^0(|\DR(\G)|) \longrightarrow \DR(\G).$$
Now we invoke Proposition \ref{p9}, which shows that $H^0(|\DR(\G)|)$ is also isomorphic
to the kernel of the morphism $\OO_X \longrightarrow \Omega^1_Y/D$, where $D$ is the differential ideal
defining the singular foliation $\tau_0(\G)$. In particular, the natural morphism $f^{-1}(\OO_S) \longrightarrow \OO_Y$ factors through
this kernel, and thus produces a morphism of sheaves of algebras
$f^{-1}(\OO_S) \longrightarrow H^0(|\DR(\G)|)$. By composition with the morphism $u$ above, we get 
the required morphism of graded mixed cdga's $v : f^{-1}(\OO_S) \longrightarrow \DR(\G)$. As explained
before this defines a natural morphism
$$\bar{v}: \DR(Y/S) \longrightarrow \DR(\G).$$
To finish the proof of the corollary, we need to show that $\bar{v}$
is, in fact, an equivalence of derived foliations. For this, it is enough to show that 
it induces an equivalence on the corresponding cotangent complexes
$$\LL_{Y/S} \longrightarrow \LL_{\G}.$$
Note that $f : Y \longrightarrow S$ being generically smooth on $Y$, the natural projections
$$\LL_{Y/S} \longrightarrow \Omega_{Y/S}^1 \qquad \LL_\G \longrightarrow \Omega_Y^1/D$$
are quasi-isomorphisms. Finally, by construction, the induced morphism $\LL_{Y/S} \longrightarrow \LL_{\G}$
is the isomorphism $\Omega^1_{Y/S} \simeq \Omega_Y^1/D$ coming from the fact that $f$ integrates the singular $\tau_0(\G)$.
\hfill $\Box$ \\

\section{Derived categories of algebraic foliations}

In this section we define and study the derived category of \emph{crystals} over a derived foliation, both in
the algebraic and in the analytic setting. In this paper we will consider only \emph{perfect crystals} on smooth varieties. The 
study of more general derived categories of derived foliations will appear elsewhere.

\subsection{Crystals along a derived foliation}\label{crys}

We let $\F \in \Fol(X)$ be a derived foliation on a smooth algebraic variety $X$
and $\DR(\F)$ the corresponding graded mixed cdga. We consider $\DR(\F)-\medg$ the 
$\s$-category of graded mixed $\DR(\F)$-dg-modules. 

\begin{df}\label{d7}
With the notations above,  a \emph{perfect crystal over $\F$}
is a graded mixed $\DR(\F)$-dg-module $E$ satisfying the following two
conditions.

\begin{itemize}
\item The dg-module $E(0)$ is perfect over $\OO_X\simeq \DR(\F)(0)$.

\item The natural morphism
$$E(0)\otimes_{\OO_X}\DR(\F) \longrightarrow E$$
is a quasi-isomorphism of graded $\DR(\F)$-dg-modules.

\end{itemize}
For a perfect crystal $E$ over $\F$, the perfect complex $E(0)$ on $X$, will be referred to as the 
\emph{underlying perfect complex of $E$}. A perfect crystal $E$ over $\F$ will be simply called a \emph{vector bundle over} $\F$ if its underlying perfect complex is quasi-isomorphic to a vector bundle on $X$ sitting in degree $0$.
The $\s$-category of perfect crystals over $\F$ is 
the full sub-$\s$-category $\Parf(\F)$ of $\DR(\F)-\medg$
consisting of perfect crystals. The full sub-$\s$-category of $\Parf(\F)$ consisting of vector bundles over $\F$ will be denoted by $\Vect(\F)$.
\end{df}

Since we will only be considering perfect crystals we will often omit the adjective perfect when
speaking about objects in $\Parf(\F)$.\\

The $\s$-category $\Parf(\F)$ is obviously functorial in $\F$ in the following sense.
Let $f : X \longrightarrow Y$ be a morphism of smooth algebraic varieties. Let 
$\F' \in \Fol(Y)$ and $\F\in \Fol(X)$ be derived foliations, and
$u : \F \longrightarrow f^*(\F') $ a morphism in $\Fol(X)$. Then, there is a base change
$\s$-functor
$$\Parf(\F') \longrightarrow \Parf(\F).$$
On affine derived foliation this base change functor is simply induced by 
the usual base change $\DR(\F) \otimes_{\DR(\F')}(-)$ on graded mixed dg-modules. \\

Before proceeding to de Rham cohomology of crystals, let us give two specific examples of 
crystals, relating this notion to more standard notions of $\D$-modules
and more generally of representations of Lie algebroids. \\

\noindent \textbf{Crystals and $\D$-modules.} Let $E$ be a quasi-coherent complex
of $\D$-modules on a smooth algebraic variety $X$. We can represent
$E$ as a pair consisting of a quasi-coherent $E(0)$ complex on $X$ together with a flat connection
$\nabla : E(0) \longrightarrow E(0)\otimes_{\OO_X}\Omega^1_X.$ The de Rham 
complex of this connection produces a graded mixed structure
on $\DR(E)=E(0) \otimes_{\OO_X} \DR(X)$, making it into a graded mixed
$\DR(X)$-module. This construction defines an equivalence between 
the $\s$-category of quasi-coherent $\D$-modules on $X$ and 
the $\s$-category of graded mixed $\DR(X)$-dg-modules which are
graded-free (see \cite[Prop. 1.1]{pato}). Restricting to perfect complexes
we see that perfect crystals over the final derived foliation  on $X$
form an $\s$-category naturally equivalent to 
the $\s$-category of $\D$-modules which are perfect over $X$. \\

\noindent \textbf{Crystals over smooth derived foliations.} Let now $\F$ be a \emph{smooth} 
derived foliation over a smooth algebraic variety $X$. We have seen already (\S \ref{liealgebroid}) that 
$\F$ corresponds to a Lie algebroid $a : T \longrightarrow \T_X$. 
A representation of this Lie algebroid is, by definition, a pair consisting of a vector bundle
$V$ together with a morphism $\nabla : V \longrightarrow V\otimes_{\OO_X}T^*$
satisfying the obvious Leibniz rule, and $\nabla^2=0$. Such a representation
has a de Rham complex $\DR(V):=V\otimes_{\OO_X}\DR(X)$, on which $\nabla$ defines
a graded mixed structure. This construction produces an $\s$-functor
from the category of representations of the Lie algebroid $T$ to 
the $\s$-category $\Parf(\F)$. It is easy to show that this
$\s$-functor is fully faithful, and that its essential image is $\Vect(\F)$. \\

The derived categories of crystals over derived foliations can be used in order
to define \emph{de Rham cohomology} of derived foliations with coefficients in a crystal. This de Rham cohomology is
usually referred to as \emph{foliated} or \emph{leafwise} cohomology in the  setting
of classical foliations or Lie algebroids. 

Let $\F$ be a derived foliation over a smooth algebraic variety $X$. The $\s$-category 
$\Parf(\F)$ can be identified with a dg-category (or, equivalently, has the structure of a
$\C$-linear $\s$-category), and we will simply denote by $Hom(E,E')$ the complex of $\C$-vectors spaces
of Hom's from $E$ to $E'$ in this dg-category structure. Notice that the projection $\DR(\F) \rightarrow \OO_X$ on the weight $0$ part, defines a structure of a graded mixed $\DR(X)$-module on 
$\OO_X$ (concentrated in weight $0$ and degree $0$).
For an arbitrary perfect crystal $E$ over $\F$, we consider the complex of morphisms
$Hom(\OO_X,E)$, from $\OO_X$ to $E$. We denote this complex of $\C$-vector spaces by
$$\HH_{DR}(\F;E):=Hom(\OO_X,E)$$
and call it \emph{de Rham cohomology of $\F$ with coefficients in $E$}. 
Note that the dg-category $\Parf(\F)$ is also endowed with a natural \emph{closed symmetric monoidal structure} (equivalently, has a structure of $\C$-linear closed symmetric monoidal $\s$-category) induced
by tensor products of $\OO_X$-modules.

\begin{rmk}\emph{Though we will not 
need this in this paper, we mention that the symmetric monoidal structure on $\Parf(\F)$
implies existence natural multiplicative structure morphisms
$$\HH_{DR}(\F;E) \otimes_{\C} \HH_{DR}(\F;E') \longrightarrow 
\HH_{DR}(\F,E\otimes_{\OO_X}E').$$
In particular, $\HH_{DR}(\F; \OO_X)$ is a commutative dg-algebra, and $\HH_{DR}(\F;E)$ 
is a dg-module over $\HH_{DR}(\F;\OO_X)$, for any $E\in \Parf(\F)$.}
\end{rmk}

The complex $\HH_{DR}(\F,E)$ can also be described as the hypercohomology
of $X$ with coefficients in an explicit complex of sheaves. We consider
$E$ as a sheaf of graded mixed $\DR(X)$-modules. We can apply the realization
$\s$-functor $|-| : \medg \longrightarrow \dg$, and thus get  
a complex of sheaves $|E|$ of $\C$-vector spaces on $X$. We then have a natural 
quasi-isomorphism
$$\HH_{DR}(\F,E) \simeq \HH(X,|E|).$$
Note that $|E|$ is explicitly given by 
the complex of sheaves $\prod_{i\geq 0}E(i)[-2i]$ endowed with 
its total differential (sum of the cohomological differential and the
mixed structure). As $E(i)$ is naturaly equivalent to $E(0)\otimes_{\OO_X}\wedge^i \LL_\F[i]$, 
$\HH_{DR}(\F,E)$ may be considered  as a version of the (derived) 
\emph{de Rham complex of $E$ along the foliation} $\F$.\\

Proposition \ref{p4} has the following extension, stating that perfect
cyrstals on quasi-smooth and rigid derived foliations are always, locally, 
pull-backs of perfect crystals on a smooth and rigid foliation on a formal scheme.

\begin{prop}\label{p7}
Let $\F \in \Fol(X)$ be a rigid and quasi-smooth derived foliation on 
a smooth algebraic variety $X$, and $E \in \Parf(\F)$ be a perfect crystal
on $\F$. Then, Zariski locally on $X$, there exists
a smooth and rigid foliation $\widehat{\F}$ on $X\times \widehat{\mathbb{A}}^n$, 
and a perfect crystal $\widehat{E} \in \Parf(\widehat{\F})$, such that 
$$e^*(\widehat{E})\simeq E$$
where $e  : X \hookrightarrow X\times \widehat{\mathbb{A}}^n$
is the zero section.
\end{prop}

\noindent\textbf{Proof.} The proof is almost the same, verbatim, as that 
of Proposition \ref{p4}. The only changes consist in considering pairs of a  
graded mixed cdga's $A$ together with a graded mixed $A$-module $M$, all along the argument. 
We leave these details to the reader. \hfill $\Box$ \\

Let $D\subset \Omega^1_X$ be a (underived) singular foliation on $X$. A \emph{coherent sheaf with 
flat connection along $D$} is defined to be a coherent sheaf $E$ on $X$ 
together with a $\C$-linear map
$$\nabla: E \longrightarrow E\otimes^u_{\OO_X}(\Omega_X^1/D)$$
satisfying the usual Leibniz rule, and being flat (that is $\nabla^2=0$ in the obvious sense). 
In terms of underived 
graded mixed algebras, coherent sheaves with flat connection along $D$ are exactly graded mixed $\DR(D)$-modules $E$ 
such that $E(0)$ is coherent and $E$ is free on $E(0)$.  If we denote by 
$\Coh(D)$ the category of the coherent sheaves with flat connection along $D$, then, there
is a \emph{truncation} $\s$-functor
$$\tau_0: \Vect(\F) \longrightarrow \Coh(\tau_{0}(\F))$$
which sends a crystal $E$ to $E(0)$ endowed with the induced map
$$E(0) \longrightarrow E(1) \longrightarrow H^0(E(1))\simeq E(0)\otimes^u_{\OO_X} (\Omega_X^1/D).$$

\bigskip

Exactly as crystals on a derived foliation have de Rham cohomology, coherent sheaves with flat connections along 
a singular foliation $D$ have naive de Rham complexes. For such an object $(E, \nabla)\in \Coh(D)$, we form
its de Rham complex modulo $D$
$$\xymatrix{
E \ar[r]^-{\nabla} & E\otimes (\Omega_X^1/D) \ar[r] & \dots \ar[r] & E \otimes (\Omega^i_X/<D>) \ar[r] & }$$
The hypercohomology of this complex on $X$ will be denoted by 
$$\HH_{DR,naive}(D,E).$$
If one starts with a derived foliation $\F$ on $X$, and $E \in \Vect(E)$, then the truncation functor induces
a canonical projection
$$\HH_{DR}(\F,E) \longrightarrow \HH_{DR,naive}(\tau_0(\F),\tau_{0}(E)),$$
which is functorial in an obvious manner. Without further assumptions on $\F$, this
morphism is not injective nor surjective in cohomology. However, we have the following result, showing 
that the closer $\F$ is to be smooth, the closer this morphism is to an equivalence.

\begin{rmk}\label{ancoh}\emph{If $D$ is an \emph{analytic} singular foliation on a complex manifold $Y$, the definitions of $\Coh(D)$, $\HH_{DR,naive}(D,E)$, for $E \in \Coh(D)$, and of the map $\HH_{DR}(\mathcal{G},\mathcal{E}) \longrightarrow \HH_{DR,naive}(\tau_0(\mathcal{G}),\tau_{0}(\mathcal{E})),$ for $\mathcal{G}$ an analytic derived foliation on $Y$ and $\mathcal{E}\in \Parf(\mathcal{G})$, all make sense and are completely analogous to the algebraic case treated above. And, obviously, the usual analytification functor for sheaves, induces a functor $\Coh(F) \to \Coh(F^{h})$, for an algebraic singular foliation $F$ on a smooth algebraic variety $X$, where $F^{an}$ is the analytification of $F$, which is an analytic singular foliation on $X^{h}$.} 
\end{rmk}

\begin{prop}\label{p9}
Let $\F$ be a quasi-smooth, rigid and derived foliation over a smooth algebraic
variety $X$ and $E\in \Vect(\F)$. Assume that $\LL_{\F}$ is quasi-isomorphic to a vector bundle on a Zariski open in $X$ whose
complement is of codimension at least $d\geq 1$. Then, the morphism
$$\HH^i_{DR}(\F,E) \longrightarrow \HH^i_{DR,naive}(\tau_0(\F),\tau_0(E)))$$
is an isomorphism for $i<d-1$ and is surjective for $i=d-1$.
\end{prop}

\noindent\textbf{Proof.} We start by the following well known lemma. 

\begin{lem}\label{lp9}
Let $\LL:=\xymatrix{V \ar[r]^-{s} & W}$ be a complex of vector bundles on $X$ with $W$ sitting in degree $0$, and 
assume that $s$ is a monomorphism which is a sub-bundle on a Zariski open $(X-S)$ with
$S$ of codimension $\geq d$. Then, for all $p\geq 0$, the perfect complex
$\wedge^p(\LL)[p]$ is cohomologically concentrated in degrees $[\inf(d-2p,-p),-p]$. 
\end{lem}

\noindent\textit{Proof of the lemma.} The complex computing $E:=\wedge^p(\LL)$ is the perfect complex
whose term in degree $-i$ is given by $Sym^{p-i}(V)\otimes \wedge^i(W)$. It is of amplitude
contained in $[-p,0]$. By assumption its higher cohomology sheaves $H^i$ for $i< 0$ are all coherent sheaves
with supports of dimension less than $(\dim X - d)$. Its dual $E^*$ is a perfect complex
of amplitude $[0,p]$ with all its higher cohomology sheaves being coherent 
with supports of dimension less than $(\dim X - d)$. Therefore, we get the following vanishing for ext-sheaves
$$\mathcal{E}xt^{i}(H^i(E^*),\OO_X)=0 \; \forall \; i < d.$$
This implies that the cohomology sheaves of the perfect complex $E\simeq (E^*)^*$ are
concentrated in $[\inf(d-p,0),0]$ as required.
\hfill $\Box$ \\

We now consider the morphism of complexes of sheaves on $X$
$$|E| \longrightarrow \HH^i_{DR,naive}(\tau_0(\F),\tau_0(E))).$$
The fiber of this map is the realization of the
graded mixed complex $K$ whose weight $p$ piece is 
$$K(p)=\tau_{\leq -1}(\wedge^p(\LL_\F) \otimes E)[p],$$ 
the fiber of the natural morphism 
$$(\wedge^p(\LL_\F) \otimes E)[p] \longrightarrow (\Omega^p_X/<D>\otimes E)[p].$$ 
By the previous lemma, each $K(p)$ sits
in cohomological degrees $[\inf(d-2p,-p),\s)$, and thus $K(p)[-2p]$ sits in
cohomological degrees $[\inf(d,p),\s)$. Note also that $K(p)\simeq 0$ if $p\leq d$.
Therefore, its realization, given by the complex
$\prod_{p> d}K(p)[-2p]$ with the appropriate differential, lies in degrees $[d,\s)$. This
concludes the proof of the proposition. \hfill $\Box$ \\

\subsection{Analytification and nilpotent crystals}

The analytification $\s$-functor for derived foliation (see \ref{analitificazione}) has the following variant 
for crystals. Let $\F$ be a derived foliation on smooth algebraic variety $X$, 
and $\F^h$ the
corresponding analytic derived foliation on $X^h$. We define
$\Parf(\F^h)$ as the $\s$-category of sheaves of graded mixed $\DR(\F^h)$-modules
$E$ on $X^h$ satisfying the following two conditions.

\begin{itemize}
\item The dg-module $E(0)$ is perfect over $\OO_{X^h}\simeq \DR(\F^h)(0)$.

\item The natural morphism
$$E(0)\otimes_{\OO_{X^h}}\DR(\F^h) \longrightarrow E$$
is a quasi-isomorphism of graded $\DR(\F^h)$-dg-modules.

\end{itemize}

It is easy to check that the usual analytification $\s$-functor for sheaves induces an $\s$-functor
$$(-)^h : \Parf(\F) \longrightarrow \Parf(\F^h).$$

Before introducing nilpotent crystals, we note that the following GAGA result for perfect 
crystals. 

\begin{prop}\label{p5'}
Let $X$ be a smooth and proper algebraic variety, and $\F$ be a derived
foliation on $X$. Then, the analytification $\s$-functor
$$\Parf(\F) \longrightarrow \Parf(\F^h)$$
is an equivalence of $\s$-categories.
\end{prop}

\noindent\textbf{Proof.} Similar to the proof of Proposition \ref{p5} and left 
to the reader. \hfill $\Box$ \\

We now introduce the \emph{nilpotent} crystals. 

\begin{df}\label{d8}
Let $X$ be a smooth algebraic variety, $\F$ a derived foliation
on $X$, and $E \in \Parf(\F)$ be a perfect crystal over $\F$. 
We say that $E$ is \emph{nilpotent} if, locally on $X^h$, the sheaf of 
of graded mixed $\DR(\F^h)$-dg-modules corresponding to 
$E^h$, belongs to the thick triangulated subcategory generated 
by the trivial crystal $\DR(\F^h)$ (considered as a graded mixed
dg-module over itself).
\end{df}

The adjective nilpotent in Definition \ref{d8} is justified by the following 
observation. Let $X=Spec\, \C$ and $\F$ be given by a finite
dimensional Lie algebra $\frak{g}$. A crystal on $\F$ is nothing else
than a complex of linear representations $E$ of $\frak{g}$ with 
finite dimensional cohomologies. Such an object is a nilpotent 
crystal if and only if furthermore the induced representation
of $\frak{g}$ on $H^i(E)$ is \emph{nilpotent} for all $i$. This shows, in particular, that
nilpotency for crystals is a non-trivial condition.
However, we will now show that it always holds for an important class of examples.

\begin{thm}\label{t0}
Let $\F \in \Fol(X)$ be a quasi-smooth and rigid derived foliation
on a smooth algebraic variety $X$. If $Sing(\F) \subset X$ is of
codimension at least $2$, then any $E \in \Vect(\F)$ is nilpotent.
\end{thm}

\noindent\textbf{Proof.} This theorem will be in fact a consequence
of Proposition \ref{p6}. Let $E \in \Vect(\F)$, and $S(E)$ the free graded mixed $\DR(\F)$-algebra
generated by $E$, i.e. 
$S(E)=Sym(E)$, where $Sym$ is computed inside the symmetric monoidal $\s$-category $\Parf(\F)$. 
The graded mixed cdga $S(E)$ defines a derived foliation $\F_E$ over $V$, where $\pi : V 
\longrightarrow X$ is the total space 
of the vector bundle $E(0)^*$ on $X$, dual of $E(0)$.  As a graded algebra, we have
$$S(E)\simeq Sym_{\OO_X}(E(0))\otimes_{\OO_X}Sym_{\OO_X}(\LL_\F[1]),$$
so that, in particular, the cotangent complex of $\F_E$ 
is $\pi^*(\LL_\F)$, hence $\F_E$ is a quasi-smooth derived foliation on $V$. \\ 
Now, the natural morphism
$\Omega^1_{V} \longrightarrow \pi^*(\LL_\F)$ is given by the derivation 
$Sym_{\OO_X}(E(0)) \longrightarrow Sym_{\OO_X}(E(0)) \otimes_{\OO_X}\LL_F$, itself
induced by the derivation $\Omega_X^1 \rightarrow \LL_\F$ and 
by 
multiplicativity from the crystal structure 
$$E(0) \longrightarrow E(0) \otimes_{\OO_X}\LL_F$$ of $E$.
In particular, the induced morphism of coherent sheaves
$\Omega^1_{V} \longrightarrow \pi^*(\LL_\F)\simeq \LL_{\F_E}$ is surjective on $H^0$
because the composite 
$$\pi^*(\Omega_X^1) \longrightarrow \Omega^1_{V} \longrightarrow \pi^*(\LL_\F)$$
is the pull-back by $\pi$ of $\Omega_X^1 \rightarrow \LL_F$ (which is
surjective on $H^0$ by the rigidity assumption on $\F$). I.e. the derived foliation $\F_E$
on $V$ is also rigid.

The derived foliation $\F_E$ on $V$ thus satisfies both the conditions of Proposition \ref{p6},
so that its truncation $\tau_0(\F^h_E)$ can be integrated locally on $V^h$.
Let $f_1,\dots,f_k$ be holomorphic functions, defined locally around a point $x\in X^h \subset V^h$,
that integrate $\tau_0(\F^h_E)$. We fix a local trivialization of $\pi^h:V^h \to X^h$ around $x$, as $X^h \times V_x$, where $V_x$ is the fiber of $V$ at $x$, so that 
local parameters at $x$ on $V^h$ are given by $(z_1,\dots,z_n,t_1,\dots,t_r)$, where 
$z_i$ are local parameters on $X^h$ and $t_j\in V_x$ form a vector space basis.
Let us consider the Taylor series expansions
of the functions $f_i$ as 
$$f_i=\sum_{\nu=(\nu_1,\dots,\nu_r)} a_{i,\nu}t^\nu,$$
where $a_{i,\nu}$ are functions around $x$ on $X^h$. \\

By construction, the homogenous part of degree $p$ of each function $f_i$ are
sections of the bundle $Sym^{p}(V^*)^h$. Moreover, the fact that 
the $f_i$ integrate the foliation $\tau_0(\F_E)$ implies that 
these sections are in fact flat sections, i.e. they are local sections of  
the sheaf $H^0_{DR,naive}(Sym^p(E)^h)$. We consider new functions at $x \in V^h$ by
taking the weight one pieces 
$$g_{i}:=\sum_{j} a_{i,(j)}t_j,$$
wherer $(j):=(0,\dots,1,\dots,0)$ and $1\leq j\leq r$.
Each function $g_i$ can be, and will be, considered
as a germ of holomorphic section at $x$ of the vector bundle $E(0)$, the dual of $V$.
These germs are flat, and thus they collectively define
an analytic germ morphisms of crystals at $x$ 
$$\phi : \OO_{X^h}^k \longrightarrow E^h.$$
The above morphism $\phi$ is clearly an isomorphism ouside of the closed subset $Sing(\F) \subset X$,
so $\phi$ is in particular an injective morphism of sheaves.
As $Sing(\F)$ is of codimension $\geq 2$,by Hartogs theorem, we deduce that $\phi$ is in fact an 
isomorphism at $x$. 

This finishes the proof that $E$ is nilpotent as required.
\hfill $\Box$ \\

\section{The Riemann-Hilbert correspondence}

\subsection{The sheaf of flat functions}\label{sflat}

Let $\F \in \Fol(X)$ be a derived foliation over a smooth algebraic variety $X$, 
and $\F^h \in \Fol(X^h)$ its analytification. The realization $|\DR(\F^h)|$ of $\DR(\F^h)$ (see Section \ref{reminder}) defines
a sheaf of commutative dg-algebras on $X^h$.

\begin{df}\label{d9}
With the notations above, the sheaf of cdga's $|\DR(\F^h)|$ on $X^h$ is called the
\emph{sheaf of locally constant functions along $\F$}. It is denoted by 
$\OO_{\F^h}$.
\end{df}

\subsubsection{The smooth case.}\label{smoothcase} In the case of \emph{smooth} derived foliations the sheaf $\OO_{\F^h}$ is easy to understand. 
As a start, let us assume that $\F$ is furthermore rigid, so that it defines 
an actual classical holomorphic regular foliation on $X^h$. The sheaf of cdga's $\OO_{\F^h}$ is then 
concentred in degree zero and is isomorphic to the subsheaf of 
$\OO_{X^h}$ of holomorphic functions which are
locally constant along the leaves. This explains the name we gave to $\OO_{\F^h}$ in general. Locally on $X^h$, the foliation is induced by a smooth
holomorphic map $f : X^h \longrightarrow \C^n$
and the sheaf $\OO_{\F^h}$ is simply given by $f^{-1}(\OO_{\C^n}^h)$. 

When $\F$ is smooth but not necessarily rigid, it corresponds to a Lie algebroid
$T \longrightarrow \T_X$ on $X$ (see \ref{liealgebroid}). Its analytification is thus a holomorphic Lie algebroid
$T^h \longrightarrow \T_X^h$ on $X^h$. This Lie algebroid possesses a 
Chevalley-Eilenberg cohomology complex $C^*(T^h)$, which is a sheaf of $\C$-linear cdga's
on $X^h$. It is explicitely given by $Sym_{\OO_{X^h}}((T^*)^h[-1])$ endowed
with the standard Chevalley-Eilenberg differential. Then, the sheaf of cdga's
$\OO_{\F^h}$ is quasi-isomprhic to $C^*(T^h)$. Therefore, for an arbitrary non-rigid smooth derived 
foliation $\F$, $\OO_{\F^h}$   is cohomologically bounded, but not necessarily concentrated in degree $0$ anymore. 
For
instance, if the Lie algebroid $T$ is abelian then $\OO_{\F^h}$ is 
$Sym_{\OO_{X^h}}((T^*)^h[-1])$ with trivial differential.

\subsubsection{Local structure.} The local structure of the sheaf $\OO_{\F^h}$  can be understood
using Proposition \ref{p4} as follows.
Let $x\in X$ and let $A=\OO_{X,x}^h$ be the
ring of germs of holomorphic functions on $X$ at $x$. By Proposition \ref{p4}
we know that there exist a smooth and rigid foliation
$\F'$ on $B=A[[t_1,\dots,t_n]]$, i.e. $k$ linearly independent
commuting derivations $\nu_1,\dots,\nu_k$ on $A[[t_1,\dots,t_n]]$ which are
holomorphic along $X$ and formal in the variables $t_i$'s. More explicitly, the derivations $\nu_i$ are given by expressions of the form
$$\sum_{i} a_i.\frac{\partial }{\partial x_i} + \sum_j b_j.\frac{\partial }{\partial t_j}$$
where $a_i$ and $b_j$ are elements of $B$ and
the $x_i$'s are local coordinates of $X$ at $x$. The derivations $\nu_i$ define
a de Rham complex
$$\xymatrix{
B\ar[r] & \oplus_i B.\frac{\partial }{\partial x_i} 
\oplus_j B.\frac{\partial }{\partial t_j}
\ar[r] & \dots
}$$
This complex is a model for the 
stalk of $\OO_{\F^h}$ at the point $x$. \\

We can say more when the derived foliation $\F$ is moreover \emph{locally d-integrable} (Definition \ref{analocint}). 
Assume that we are given $f : X \longrightarrow Y$ a holomorphic
map between complex manifolds. As we are interested
in a local description we can assume that $X$ and $Y$ are Stein manifolds, 
and will allow ourselves to restrict to open subsets in necessary.
We assume that $\F:=f^*(\mathbf{0}_Y)$ is the induced
analytic derived foliation on $X$. By construction, the sheaf of cdga $\OO_{\F}$
is the sheaf on $X$, for the analytic topology, of relative 
derived de Rham cohomology: it associates to 
an open $U \subset X$ the cdga $\OO_{\F}(U):=|\DR(\F)|=|\DR(U/Y)|$. This
sheaf can explicitly be computed as follows. We let $\xymatrix{
X \ar[r]^-{j} & Z=X\times Y \ar[r]^-p & Y}$ be a factorization of $f$
as a closed immersion followed by a smooth morphism (where
$j$ is the graph of $f$). The sheaf 
of rings $\widehat{\OO}_{X}$ of formal functions on $Z$ along $X$, 
has a natural structure of a $\D_Z$-module on $Z$. Therefore, we can consider the
relative de Rham complex of $\widehat{\OO}_X$ over $Y$, and obtain  
a sheaf $|\DR(\widehat{\OO}_X/Y)|$ of cdga's on $Z$. This sheaf  is set-theoretically supported on $X$, 
and thus can be considered as a sheaf of cdga's over $X$.

\begin{lem}\label{l1}
There exists an equivalence of sheaves of cdga's on $X$
$$|\DR(\widehat{\OO}_X/Y)| \simeq \OO_{\F}.$$
\end{lem}

\noindent\textbf{Proof.} The cotangent complex of $\F$ sits in the following
exact triangle of perfect complexes on $X$
$$\xymatrix{
N^* \ar[r] & \Omega_X^1 \ar[r] & \LL_\F,}$$
where $N^*=f^*(\Omega_Y^1)$. 
We thus have a push-out of graded mixed cdga's over $X$
$$\xymatrix{
\DR(\F) \ar[r] & Sym_{\OO_X}(N^*[2])=:B \\
\DR(X) \ar[r] \ar[u] & \OO_X \ar[u]
}$$
As $X$ is Stein, we see that the graded mixed structure on $Sym_{\OO_X}(N^*[2])$
is automatically trivial. We consider the internal de Rham algebra
$\DR^{int}(B/\DR(\F))$. This is a cdga endowed with two 
extra gradings, and two commuting graded mixed structures. As such, it
has a total realization defined by
$$||\DR^{int}(B/\DR(\F))||:=\rch(\OO_X,\DR^{int}(B/\DR(\F)))$$
where $\OO_X$ is considered as a trivial double graded mixed complex of sheaves.
This total realization can be obtained by successive realizations of the two
graded mixed structures, and thus in two different manners depending
on the orders in which these realizations are taken. 
If we realize the first graded mixed structure, the one defined
by the graded mixed structure on $\DR(\F)$, 
we get $|\DR(\widehat{\OO}_X/Y)|$. Realizing the second one, gives back
$\DR(\F)$. We thus obtain this way the desired equivalence after taking
the second realization. \hfill $\Box$ \\

A consequence of the above lemma is the following explicit 
description of the stalks of $\OO_{\F}$ at a given 
point $x \in X$. Let $(x_1,\dots,x_k)$ and $(y_1,\dots,y_n)$
be local coordinates of $X$ at $x$ and $Y$ at $y=f(x)$. The ring
of germs of functions on $X\times Y$ at $(x,y)$ can then be identified
with $\C\{x_*,y_*\}$, the ring of germs of holomorphic functions
at $0 \in \C^{k+n}$, and the graph of $f$ defines an ideal $I \subset 
\C^{k+n}$. We thus have the relative de Rham 
complex of $\C\{x_*,y_*\}$ over $\C\{y_*\}$, and its completion along $I$
$$\xymatrix{
\widehat{\C\{x_*,y_*\}} \ar[r] & \widehat{\C\{x_*,y_*\}}^k \ar[r]
& \wedge^2(\widehat{\C\{x_*,y_*\}}^k) \ar[r] & \dots &
\wedge^k(\widehat{\C\{x_*,y_*\}}^k).}$$
The differential in this complex is induced by 
the derivative relative to $Y$, sending
a function $f$ to $\sum_{j\leq k}\frac{\partial f}{\partial x_j} . dx_j$.
The above complex is naturaly quasi-isomorphic to the stalk
of $\OO_\F$ at $x$.  \\

\subsubsection{Flat functions and singularities.} In general, the complex of sheaves of flat functions $\OO_{\F}$ is not cohomologically concentrated in degree $0$, as we will show in some specific examples below. Its non-zero cohomology sheaves contain
subtle informations about the singularities of $\F$, that can be sometimes
expressed in terms of vanishing cohomology.

Let us have a look at the specific case of codimension one
quasi-smooth derived foliations, and the local structure of the sheaf $\OO_{\F^h}$.
Assume that locally, in the analytic topology, such a derived foliation 
is induced by a holomorphic map
$$f : X \longrightarrow \C.$$
The derived foliation on $X$ is then $f^*(0)$, the pull-back of the trivial
foliation on $\C$.
If the map $f$ is constant, say $0$, then it is easy to see that 
$\OO_{\F^h})$ is the constant sheaf on $X$ with stalks $\C[[x]]$ the formal
functions on $\C$ at $0$. Indeed, in this case the cotangent
complex of $f$ splits as $\Omega_X^1\oplus \OO_X[1]$, and the
graded mixed cdga $\DR(X/\C)$ is $\DR(X)[u]$ where $u$ is a free variable
in degree $2$. The realization of this graded mixed cdga is
$|\DR(X)|[[x]]$ the formal power series in the absolute de Rham complex
of $X$. As a sheaf on $X$ this is quasi-isomorphic to $\C[[x]]$. This situation
generalizes easily to the case of an arbitrary holomorphic map $f : X \longrightarrow \C^n$
having constant rank (left to the reader).

Let us now assume that $f : X \longrightarrow \C$ is not constant. The graded
mixed cdga corresponding to $\F=f^*(0)$ is the relative de Rham algebra $\DR(X/\C)$. 
The cotangent complex of $f$ is now represented by the complex
$$\LL_{f}=\xymatrix{\OO_X \ar[r]^-{df} & \Omega_X^1,}$$
given by the differential of $f$. 
As graded cdga, $\DR(X/\C)$ is then $\DR(X)[u]$, with $u$ in degree $2$, and 
where the cohomological differential 
sends $u$ to $df$. The graded mixed structure on $\DR(X)[u]$, at least locally on $X$, 
is simply given by the derivation $\OO_X \longrightarrow \Omega_X^1 \longrightarrow \LL_f$
where the second map is the canonical map. According to this, the realization 
of $\DR(X/\C)$ is a version of the formal twisted de Rham complex of \cite{sab}.
To be more precise, we consider the sheaf of graded algebras on $X$
$$\widehat{Sym}_{\OO_X}(\OO_X \oplus \Omega_X^1[-1])\simeq \prod_{p\geq 0}\Omega_X^{*\leq p},$$
where $\Omega_X^*$ stands for $Sym_{\OO_X}(\Omega_X^1[-1])$.
We consider the differential
$t.d_{DR}+\wedge df$, 
where $t.d_{DR}$ means the de Rham differential going from the $p$-component
to the $(p+1)$-component.
This complex is naturally quasi-isomorphic to the realization of $\DR(\F)=\DR(X/\C)$
and thus is a model for the sheaf $\OO_{\F^h}$. 

Assume that $f$ has an isolated singularity at $x \in X$, and let us restrict
$X$ so that $x$ is the only critical point of $f$. Then, lemma \ref{lp9} implies that
for any $p$ the complexes of sheaves 
$$\xymatrix{
\OO_X \ar[r]^-{df} & \Omega_X^1 \ar[r]^-{df} & \dots \ar[r]^-{df} & \Omega_X^p}$$
are cohomologically concentrated in degree $p$ and thus quasi-isomorphic to 
the a $p$-shift of the sheaf $\mathcal{H}_p(f):=\Omega_X^p/df\wedge \Omega_X^{p-1}$. 
The fiber of the truncation map
$$|\DR(X/\C)| \longrightarrow (\Omega_X^*/df,d_{DR}),$$
where $(\Omega_X^*/df,d_{DR})$ is the naive relative de Rham complex, 
is thus  the realization of the graded mixed complex
whose component of weight $p$ are zero if $p\leq Dim X$ and
$\mathcal{H}_{Dim X}[2p-Dim X]$ if $p>n$. Therefore, this realization 
is quasi-isomorphic to $\prod_{p>Dim X}\mathcal{H}_{dim X}(f)[-dim X]$.
As obviously $H^{Dim X}(\Omega_X^*/df,d_{DR}))\simeq 0$,
we deduce that we have
$$H^p(\OO_{\F^h}) \simeq H^p(\Omega_X^*/df,d_{DR})$$
for $p<Dim X-1$, and we have a long exact sequence
$$\xymatrix{
0 \ar[r] & H^{Dim X-1}(\OO_{\F^h}) \ar[r] &  H^{Dim X-1}(\Omega_X^*/df,d_{DR})
\ar[r] &}$$
$$\xymatrix{\ar[r] &
\prod_{p> Dim X-1}\mathcal{H}_{dim X}(f) \ar[r] & 
H^p(\OO_{\F^h}) \ar[r] & 0.}$$
In particular, we see that 
the top cohomology sheaf of $\OO_{\F^h}$ is supported at $x$ and contains and infinite
number of copies of
$J(f)$, the Jacobian ring of $f$ at $x$. 

In the general situation, it is surely possible to relate the cohomology sheaves
of $\OO_{\F^h}$ to some twisted de Rham complex as described in \cite{sab}, and thus
to vanishing cohomology. Note also that the complex $\OO_{\F^h}$ appears
in a disguised form, under the name \emph{Koszul-de Rham algebra}, in \cite{sait}, where it 
is explicitely related to the naive relative de Rham complex and the singularities of the
map $f$. \\

\subsection{The Riemann-Hilbert correspondence}

Let $X$ be a smooth algebraic variety and $\F \in \Fol(X)$ be a 
derived quasi-smooth foliation on $X$. We first construct the \emph{Riemann-Hilbert $\s$-functor} for $\F$
$$RH : \Parf(\F) \longrightarrow \OO_{\F^h}-\dg.$$
It is defined as the composite 
$$\xymatrix{\Parf(\F) \ar[r]^-{(-)^h} & \Parf(\F^h) \ar[rrr]^-{\rch(\OO_{X^h},-)} &&& \OO_{\F^h}-\dg  }$$
where $\rch(\OO_{X^h},-)$  sends a crystal $E$ over $\F^h$ to the complex 
$\rch(\OO_{X^h},E)$ of sheaves of $\C$-vector spaces, endowed with its natural module structure over
$$\rch(\OO_{X^h},\OO_{X^h})\simeq |\DR(\F^h)|=\OO_{\F^h}.$$ We let 
$\Parf^{nil}(\F)$ be the full sub-$\s$-category of 
$\Parf(\F)$ consisting of nilpotent crystals in the sense of Definition \ref{d8}.
By definition of nilpotency for crystals, the $\s$-functor $RH$ defined above 
restricts to 
$$RH : \Parf^{nil}(\F) \longrightarrow \Parf(\OO_{\F^h}),$$
where $\Parf(\OO_{\F^h})$ is the $\s$-category of sheaves of perfect
$\OO_{\F^h}$-dg-modules on $X^h$. 

\begin{thm}\label{t1}
Let $\F \in\Fol(X)$ be a quasi-smooth derived foliation on a smooth 
algebraic variety $X$. If $X$ is proper, then the $\s$-functor
$$RH : \Parf^{nil}(\F) \longrightarrow \Parf(\OO_{\F^h})$$
is an equivalence.
\end{thm}

\noindent\textbf{Proof.} By GAGA (see Proposition \ref{p5'}) we know that the 
analytification $\s$-functor induces an equivalence
$\Parf^{nil}(\F) \simeq \Parf^{nil}(\F^h)$. Now, both $\s$-categories
$\Parf^{nil}(\F^h)$ and $\Parf(\OO_{\F^h})$ are global sections of natural
stacks of triangulated $\s$-categories on $X^h$.
Let us denote these 
stacks by $\underline{\Parf}^{nil}(\F^h)$ and $\underline{\Parf}(\OO_{\F^h})$.
Moreover, the $RH$ $\s$-functor is itself obtained by taking global sections of
a morphism of stacks $\underline{RH} : \underline{\Parf}^{nil}(\F^h)
\longrightarrow \underline{\Parf}(\OO_{\F^h})$.
By definition of nilpotency, $\OO_{X^h}$ locally generates (in the sense of triangulated $\s$-categories) $\underline{\Parf}^{nil}(\F^h)$, while $\OO_{\F^h}$ locally generates $\Parf(\OO_{\F^h})$, and by definition $RH(\OO_{X^h})\simeq \OO_{\F^h}$. We then conclude that $\underline{RH}$ is locally fully faithful and locally
essentially surjective. It is therefore a local equivalence of stacks, and thus
a global equivalence, i.e. $RH$ is an equivalence. \hfill $\Box$\\

Combining Theorem \ref{t1} with Theorem \ref{t0}, we get 
our main result. For $\F \in\Fol(X)$ be a quasi-smooth and rigid derived foliation on a smooth and proper
algebraic variety $X$, we denote by $\Parf^{v}(\F) \subset \Parf(\F)$
the full triangulated sub-$\s$-category generated by objects in $\Vect(\F)$.

\begin{cor}\label{ct1}
Let $\F \in\Fol(X)$ be a quasi-smooth and rigid derived foliation on a smooth and proper
algebraic variety $X$. If the singular locus $Sing(\F)$ has codimension $\geq 2$, then 
the Riemann-Hilbert correspondence
$$RH : \Parf^{v}(\F) \longrightarrow \Parf^v(\OO_{\F^h})$$
is an equivalence of $\s$-categories.
\end{cor}

\subsection{Examples}

\subsubsection{Smooth and rigid foliations.}\label{deligne} Let $X$ be a smooth and proper algebraic variety and 
$\F$ be \emph{smooth and rigid} derived foliations on $X$. We represent 
$\F$ by a differential ideal $D \subset \Omega_X^1$ which is furthermore a subbundle. 
The category $\Vect(\F)$, of vector bundle crystals along $\F$ is equivalent to the category of
vector bundles $V$ on $X$ endowed with a partial flat connection
$$\nabla : V \longrightarrow V\otimes \Omega_X^1/D$$
satisfying the standard properties. As already observed (see \ref{smoothcase}), the sheaf $\OO_{\F^h}$ is then a genuine sheaf of commutative algebras:
it is the subsheaf of $\OO_{X^h}$, of holomorphic functions on $X$, which are locally constant
along the leaves, i.e. local functions $f$ on $X$ such that $df \in D \subset \Omega_X^1$. Locally on $X^h$, 
the foliation is given by a smooth holomorphic map $X \longrightarrow \C^p$, and the
sheaf $\OO_{\F^h}$ can then be described as $f^{-1}(\OO_{\C^p})$, the sheafy inverse image of the sheaf of holomorphic
functions on $\C^p$. 

The Riemann-Hilbert correspondence of Theorem \ref{t1} and Corollary \ref{ct1} implies the existence of an equivalence of categories
$$\Vect(\F) \simeq \Vect(\OO_{\F^h})$$
where $\Vect(\OO_{\F^h})$ is the category of sheaves of $\OO_{\F^h}$-modules on $X^h$ which are locally
free of finite rank. When $\F$ is furthermore induced by a smooth and proper morphism
$f : X \longrightarrow Y$, i.e. $f^*(\mathbf{0}_Y)\simeq \F$, this equivalence is essentially the relative
Riemann-Hilbert correspondence of \cite{del}. The extension from vector bundles to perfect complexes
essentially states that this equivalence is also compatible with cohomology, and also induces
a quasi-isomorphism
$$\HH_{DR}(\F;V) \simeq \HH_{B}(X^h,RH(V)),$$
between the algebraic de Rham cohomology of $V$ over $X$ to the
(Betti) cohomology of $X^h$ with coefficients in the sheaf of $\OO_{\F^h}$-modules $RH(V)$.

\subsubsection{Lie algebroids.} Let $X$ be a smooth and proper algebraic variety and 
$\F$ be smooth derived foliation on $X$. We have seen that such a derived foliation corresponds to a classical Lie
algebroid with anchor map $a : T \longrightarrow \T_X$ (see \S \ref{liealgebroid}). The category 
$\Vect(\F)$ can then be described as representations of $T$ in vector bundles, i.e. 
vector bundles $V$ together with a connection along $\F$
$$\nabla : V \longrightarrow V\otimes T^*$$
satisfying the Leibniz and flatness $\nabla^2=0$. In this case, the sheaf $\OO_{\F^h}$ is the sheaf of cdga's
on $X^h$
$$\xymatrix{
\OO_{X^h} \ar[r] & (T^*)^h \ar[r] & \wedge^2(T^*)^h \ar[r] & \dots \ar[r] & \wedge^{d}(T^*)^h}$$ 
where $d$ is the rank of $V$. This is a sheaf of cdga's not concentrated in degree 
$0$ in general. If the anchor map $a$ turns out to be a subbundle on a Zariski open $U \subset X$, 
then the higher cohomology sheaves have supports on $X-U$. 

We denote by $\mathcal{K}$ the kernel of the anchor map $a : T \longrightarrow \T_X$, considered
as a sheaf on the big \'etale site of $X$. As such it is representable by 
a Lie algebra scheme $K \longrightarrow X$, whose fiber $K_x$ at a point $x \in X$ is the
kernel of the map $a_x : T_x \longrightarrow \T_{X,x}$, which is a $\C$-linear Lie
algebra. For an object $V \in \Vect(\F)$, and a point $x\in X$, the Lie algebra $K_x$ acts
naturally on the fiber $V_x$. If $V$ is nilpotent in the sense of 
Definition \ref{d8}, then for all points $x \in X$ the $K_x$-module
$V_x$ is a nilpotent representation of the Lie algebra $K_x$. 

The Riemann-Hilbert correspondence of Theorem \ref{t1} induces an equivalence of categories
$$\Vect^{nil}({\F}) \simeq \Vect(\OO_{\F^h}).$$
The extension to perfect complexes again essentially states that this equivalence is also
compatible with the natural cohomologies on both sides, i.e. algebraic de Rham cohomology
and Betti cohomology.

\subsubsection{Log vector fields.}\label{log} Let $X$ be a smooth and proper algebraic variety and 
$D \subset X$ be a divisor with simple normal crossings. We let $\T_X(log(D)) \subset \T_X$ be the
subsheaf of vector fields that stabilize $D$, i.e. derivations preserving the ideal defining $D$. 
This defines a Lie algebroid on $X$, and thus a smooth derived foliation $\F$, 
to which we can apply our Riemann-Hilbert correspondence (Theorem \ref{t1} and Corollary \ref{ct1}). 
Let us unravel the two sides of the Riemann-Hilbert correspondence in this case.
Note first that $\Vect(\F)$ consists of pairs consisting of a vector bunlde on $X$ together with a
logarithmic flat connections \emph{along $D$}
$$\nabla : V \longrightarrow V\otimes \Omega^1_X(log(D)).$$
The structure sheaf $\OO_{\F^h}$ is then the sheaf of cdga's on $X^h$ 
given by the \emph{logarithmic de Rham complex} 
$$\xymatrix{
\OO_{X^h} \ar[r] &  \Omega^1_{X^h}(log(D)) \ar[r] &  \Omega^2_{X^h}(log(D)) \ar[r]
& \dots \ar[r] & \Omega_{X^h}^d(log(D)).}$$
By Grothendieck's log-de Rham theorem, this is a resolution of the 
sheaf of cdga's $j_*(\C)$, where $j : (X-D) \longrightarrow X$
and $\C$ is the constant sheaf (with stalks $\C$) on $(X-D)$. Note that $j_*(\C)$ is not 
concentrated in degree $0$, as its fiber at a point $x\in D$ 
is the cohomology algebra of an $m$-dimensional torus $H^*((S^1)^m,\C)$, 
if the local equation of $D$ at $m$ is of the form $x_1.\dots.x_m=0$
(for $(x_1, \cdots, x_{d}) $ a system of local parameters at $x$ in $X$, and $m\leq d$).

The Lie algebroid $\T_X(log(D))$ has isotropy along $D$. For a point $x \in D$, 
in a neigborhood of which $D$ has equation $x_1.\dots.x_m=0$, 
the kernel of the anchor map $a: \T_X(log(D))\to \T_X $ at $x$ is an 
abelian Lie algebra of dimension $m$. Generators of this Lie algebra are given 
by the local vector fields $x_i.\frac{\partial}{\partial x_i}$ for $1\leq i\leq m$, 
which are local sections of $\T_X(log(D))$ that, once evaluated at $x$, provide a basis
for the Lie algebra $Ker(a_x)$. 
As we have already remarked, if $V \in \Vect(\F)$ is nilpotent, then the actions
of the Lie algebras $Ker(a_x)$ on $V_x$ must be nilpotent. In our case, 
these actions are given by the residues of the connection $\nabla$ along the components
of $D$, and thus when $V$ is nilpotent these residues must be nilpotent too. The converse is true: 
$V\in \Vect(\F)$ is nilpotent if and only if it has nilpotent residues. This follows from the local
analytic structure of flat connections with logarithmic poles along $D$, and the fact that they
correspond to local systems on $X-D$ with unipotent local monodromies around $D$ (see e.g. \cite{del}).
The Riemann-Hilbert correspondence of Theorem \ref{t1} then becomes the following equivalence of categories
$$RH : \Vect^{nil}(\F) \simeq \Vect(j_*(\C)),$$
where $\Vect^{nil}(\F)$ can be identified with vector bundles on $X$ endowed
with flat connections with logarithmic poles along $D$. Its extension to perfect complexes
$$RH : \Parf^{nil}(\F) \simeq \Parf(j_*(\C))$$
implies that the previous equivalence $RH : \Vect^{nil}(\F) \simeq \Vect(j_*(\C))$ 
is also compatible with cohomology theories on both sides.
Note here that 
$j_*(\C)$ can also be described as $\pi_*(\C_{X^h})$, where
$\pi : X^{log} \longrightarrow X^h$ is the logarithmic homotopy type of the pair
$(X,D)$ in the sense of \cite{kana}, and that $\C$ generates, in the sense of triangulated
categories, the categrory of unipotent local systems. Therefore, perfect $\OO_{\F^h}$-dg-modules
can also be understood as perfect complexes of $\C$-modules on $X^{log}$ which are
relatively unipotent over $X^h$, i.e. those which become unipotent along the fibers of $\pi$.
These statements have a straightforward generalization
to the case where $(X,D)$ is replaced by a more general log-structure on $X$.
The equivalences above recovers the logarithmic Riemann-Hilbert
correspondence of \cite[Thm. 0.5]{kana}. 

\subsubsection{The RH correspondence along a non-smooth morphism.} Let $f : X \longrightarrow Y$ be a 
morphism between smooth and proper algebraic varieties, and $\F=f^*(0_Y)\in \Fol(X)$ be the derived
foliation induced by $f$ (recall that $f^*$ denotes here the pull-back functor on foliations, and 
$0_Y$ is the final foliation on $Y$). We associate to $f$ the morphism between de Rham shapes (first introduced in \cite{sim})
$f_{DR}: X_{DR} \longrightarrow Y_{DR}$, where for a scheme $Z$ the functor $Z_{DR}$ sends
an algebra $A$ to $Z(A_{red})$. The \emph{relative de Rham shape} is defined  
$$(X/Y)_{DR}:=X_{DR}\times_{Y_{DR}}Y.$$
Quasi-coherent sheaves on $(X/Y)_{DR}$ are by definition relative crystals on $X$ over $Y$. Note that these
are also the quasi-coherent sheaves on the \emph{relative infinitesimal site} $(X/Y)_{\textrm{inf}}$ of $X$ over $Y$ (see \cite{gro}).
This site has objects commutative diagrams of the form
$$\xymatrix{
S_{red} \ar[r] \ar[d] & X \ar[d] \\
S \ar[r] & Y,}$$
where $S_{red} \longrightarrow X$ is a Zariski open. The topology is itself defined in the natural way. 
It comes equiped with a structure sheaf $\OO_{X/Y,\textrm{inf}}$ sending a diagram as above to 
$\OO(S)$. It is then possible to prove that there exists a natural equivalence of 
$\s$-categories
$$\Parf(\F) \simeq \Parf(\OO_{X/Y,\textrm{inf}}).$$
Therefore, perfect crystals along the foliation $\F$ are nothing else than perfect 
complexs on the relative infinitesimal site of $X$ over $Y$.

In this situation, the sheaf $\OO_{\F^h}$ on $X^h$ is the relative analytic derived de Rham complex
of $X^h$ over $Y^h$. As explained earlier, it has a projection onto the naive
de Rham complex $(\Omega_{X^h/Y^h}^*[-*],dR)$, which is the de Rham complex
of relative holomorphic differential forms
$$\xymatrix{
\OO_{X^h} \ar[r] & \Omega_{X^h/Y^h }^1 \ar[r] & \dots \ar[r] & \Omega_{X^h/Y^h}^n.}$$
It is important to remark that the projection
$$\OO_{\F^h} \longrightarrow (\Omega_{X^h/Y^h}^*[-*],dR)$$
is far from being an equivalence in general. For example, if $f: X \to Y$ is a closed
immersion the right hand side is just the sheaf $\OO_{X^h}$ whereas
the left hand side is $\widehat{\OO}_{X^h}$, the structure sheaf of the formal
completion of $X$ inside $Y$. Proposition \ref{p9} implies that this morphism is
an isomorphism on cohomology groups in degree less than $d$ if the map $f$ is
smooth outside of codimension $(d+1)$ closed subset. 

Let us assume now that $f : X \longrightarrow Y$ is \emph{flat with reduced fibers and of strictly positive
relative dimension}. By generic smoothness $f$ is smooth on $X-S$ where $S$ is a closed
subset of codimension $d>1$. By Proposition \ref{p9} we have that 
$$H^0(\OO_{\F^h}) \simeq H^0((\Omega_{X^h/Y^h}^*[-*],dR)).$$
The right hand side is the subsheaf of functions on $X^h$ which are locally constant 
along the fibers of $f$, i.e. those which are pull-backs of functions on $Y^h$. In other
words, we have an isomorphism of sheaves of rings
$$H^0(\OO_{\F^h}) \simeq f^{-1}(\OO_{Y^h}).$$
Therefore, $\OO_{\F^h}$ has a canonical structure of a sheaf of 
$f^{-1}(\OO_{Y^h})$-cdga's. Therefore, for any 
point $y \in Y$, we can consider the following sheaf of cdga's on $X^h$
$$\OO_{\F^h} \otimes_{f^{-1}(\OO_{Y^h})}\C_y,$$
where the map $f^{-1}(\OO_{Y^h}) \longrightarrow \C_y=\C$ is given by evaluation at $y$. The resulting 
sheaf is the derived analytic de Rham cohomology of the fiber $f^{-1}(y)$, and thus
it is the constant sheaf $\C$ on that fiber. We deduce that
$$\OO_{\F^h} \otimes_{f^{-1}(\OO_{Y^h})}\C_y \simeq (i_y)_*(\C)$$
where $i_y : f^{-1}(y) \hookrightarrow X$ is the inclusion of the fiber at $y$.

The previous discussion prompts the following interpretation. For any $E \in \Vect(\OO_{\F^h})$, we define its base changes
$E_y := E \otimes_{f^{-1}(\OO_{Y^h})}\C_y$ which are \emph{genuine} local systems 
on $f^{-1}(y)$. Therefore, the notion of perfect complexes 
of $\OO_{\F^h}$-dg-modules should be understood as a notion of \emph{analytic families
of perfect complexes of local systems along the fibers}. With this interpretation, the Riemann-Hilbert correspondence
$$RH : \Parf(\F) \simeq \Parf(\OO_{\F^h})$$ of Theorem \ref{t1} 
should be actually understood as an equivalence between 
\emph{algebraic families of crystals along the fibers of} $f$, and \emph{analytic families of
perfect local systems along the fibers.}

\bibliographystyle{alpha}
\bibliography{Fol-Bib.bib}

\end{document}